\documentclass[11pt]{amsart}

\usepackage{amsmath}
\usepackage{amsfonts}
\usepackage{amssymb}
\usepackage{amsthm}
\usepackage{mathtools}
\usepackage{latexsym}

\usepackage{enumerate}
\usepackage{cancel}
\usepackage{cases}
\usepackage{empheq}
\usepackage{multicol}

  \usepackage{caption}
  \usepackage{subcaption}
  \usepackage{graphics}
  \usepackage{epsfig}
\usepackage{graphicx}

\usepackage{float}

\usepackage{color}

\usepackage[backgroundcolor=gray!30,linecolor=black]{todonotes}
   
\usepackage{mathrsfs}
\usepackage{fontenc}
\usepackage{inputenc}

\usepackage{verbatim}

\theoremstyle{plain}
\newtheorem{theorem}{Theorem}[section]

\newtheorem{lemma}[theorem]{Lemma}
\newtheorem{corollary}[theorem]{Corollary}

\theoremstyle{definition}
\newtheorem{definition}[theorem]{Definition}

\theoremstyle{remark}
\newtheorem{remark}[theorem]{Remark}

\numberwithin{equation}{section}
\numberwithin{figure}{section}

\usepackage[square,comma,numbers,sort]{natbib}
\usepackage[colorlinks=true, pdfborder={ 0 0 0}]{hyperref}
\hypersetup{urlcolor=blue, citecolor=red}
\usepackage{url}

\newcommand{\field}[1]{\mathbb{#1}}

\newcommand{\nR}{\field{R}}

\newcommand{\nT}{\mathbb T}

\newcommand{\ip}[2]{\left<#1,#2\right>}

\newcommand{\norm}[1]{\left\|#1\right\|}

\newcommand{\normLp}[2]{\left\|#2 \right\|_{L^{#1}}}

\newcommand{\normL}[1]{\|#1\|_{L^2}}
\newcommand{\normH}[1]{\|#1\|_{H^1}}
\newcommand{\fournorm}[1]{\|#1\|_{L^4}}
\newcommand{\infnorm}[1]{\|#1\|_{L^{\infty}}}
\newcommand{\Htwonorm}[1]{\|#1\|_{H^{2}}}

\usepackage{url}

\usepackage{placeins}

\newcounter{my_counter}
\title[3D VVV-MHD]{On a Partial Voigt Regularization of the 3D Magnetohydrodynamic Equations in Velocity-Vorticity Form  
}
\date{\today}

\author{Adam Larios}
\address[Adam Larios]{Department of Mathematics, 
                University of Nebraska--Lincoln,
        Lincoln, NE 68588-0130, USA}
\email[Adam Larios]{alarios@unl.edu}

\author{Yuan Pei}
\address[Yuan Pei]{Department of Mathematics, 
                Western Washington University
        Bellingham, WA 98225, USA}
\email[Yuan Pei]{peiy@wwu.edu}


\begin{document}
\begin{abstract}
The Velocity-Vorticity (VV) formulation of the incompressible Navier-Stokes equations has become popular in recent years, especially in numerical studies, due to its structural advantages.  Recently, with L. Rebholz, we introduced a Voigt regularization to the momentum equation in this formulation, establishing global well-posedness of the regularized system in 3D, along with convergence results and a blow-up criterion.  In the present work, we extend these ideas to the 3D magnetohydrodynamics (MHD) equations.  While it may seem that a ``VV-type'' split on the magnetic equation is required, we show that no such modification is necessary, and global well-posedness holds with a Voigt regularization only on the momentum equation, preserving the structure of both the vorticity and magnetic equations.  We also prove that the regularized system converges to the original system, up to a possible blow-up time, and we establish a blow-up criterion for solutions to the original 3D MHD system.
\end{abstract}

\maketitle
\thispagestyle{empty}

\section{Introduction}\label{secInt}
\noindent
The interaction between fluid motion and magnetic fields plays a central role in systems ranging from liquid-metal flows and fusion plasmas to astrophysical dynamics.  These phenomema are widely understood be governed by the magnetohydrodynamic (MHD) equations.  However, as we discuss below, these equations are notoriously difficult to solve accurately and efficiently due to nonlinear interactions that create a feedback between very small and very large scales, necessitating direct numerical simulations to run at very high resolutions, which is both costly and time consuming.  Moreover, it is currently unknown whether these nonlinear interactions in solutions to the MHD equations can create arbitrarily small length scales to the extent that a shock develops. If so, in experimental designs that do not have shock development, it may be of use to consider modifications of the MHD system that prevent shocks from occurring; however, this is difficult in practice, as many standard modifications such as hyper-diffusive models or filtering by convolution give rise to non-physical steady-states of the equation, and additionally often require non-physical boundary conditions to be well-posed.  

In recent years, the \textit{Voigt\footnote{As discussed in \cite{Brzezniak_Larios_Safarik_2025_JMFM}, Voigt is sometimes incorrectly spelled ``Voight'' in the literature.} regularization}\footnote{The Navier-Stokes-Voigt equations were first proposed in 1966 by O. Ladyzhenskaya in \cite{Ladyzhenskaya_1966_Voigt_talk} and later studied by her student A.~P.~Oskolkov in \cite{Oskolkov_1973,Oskolkov_1976,Oskolkov_1982} as a model for Kelvin-Voigt fluids.  They were later recognized as a regularization for the Navier-Stokes and Euler equations in \cite{Cao_Lunasin_Titi_2006}, where the Euler-Voigt equations were first considered.  Oskolkov introduced the name ``Voigt'' for this model, but W. Voigt (1950--1919) was not involved in the creation of the equation bearing his name.  We point this out as the misattribution of the Navier-Stokes-Voigt equations to W. Voigt is quite common, and appears, e.g., as recently as 2024 in  \cite{Yang_Huang_He_2024_IJNAM_Voigt_triple}.} 
has become increasingly popular as a modification to fluid models that does not suffer from the aforementioned difficulties.  In particular, it has exactly the same steady states as the original system, and moreover is globally well-posed in the presence of physical boundary conditions (in particular, one can prove that no shock develops). Moreover, the Voigt regularization even allows for global well-posedness in the infinite Reynold's number case \cite{Cao_Lunasin_Titi_2006}. In addition, it is fairly straight-forward to implement in standard fluid solvers 
\cite{DiMolfetta_Krstlulovic_Brachet_2015,
Kuberry_Larios_Rebholz_Wilson_2012,
Layton_Rebholz_2013_Voigt}, 
requiring one additional elliptic solve per evaluation of the right-hand side.  The Voigt regularization has been studied in a wide variety of analytical contexts, too numerous to discuss in detail here.  See, e.g., 
\cite{Berselli_Bisconti_2012,
Berselli_Kim_Rebholz_2016,
Berselli_Kim_Spirito_2016_DCDSB,
Berselli_Spirito_2017,
Bisconti_Catania_2021, 
Bohm_1992,
Borges_Ramos_2013,
Boutros_Liu_Thomas_Titi_2026_sea_ice,
CotiZelati_Gal_2015_Voigt,
Cuff_Dunca_Manica_Rebholz_2015,
deOliveira_Borges_Khompysh_Shakir_2025_JMP_Voigt_compressible,
DiPlinio_Giorgini_Pata_Temam_2018,
Ebrahimi_Holst_Lunasin_2012,
Gal_Medjo_2013_MMAS,
Gao_Sun_2012,
Garcia_Luengo_Julia_Read_2012,
Ilyin_Kalantarov_Zelik_2026_Voigt_attractor_dim,
Ilyin_Zelik_2025_Voigt_attractors,
Kalantarov_Levant_Titi_2009,
Kalantarov_Titi_2009,
Khouider_Titi_2008,
Krasnoschok_Pata_Siryk_Nataliya_2020_PhysD_subdiffusive,
Larios_Lunasin_Titi_2015,
Larios_Pei_Rebholz_2018,
Larios_Petersen_Titi_Wingate_2015,
Levant_Ramos_Titi_2009,
Li_Qin_2013,
Niche_2016_JDE,
Pei_2021,
Ramos_Titi_2010,
Rong_Fiodilino_Shi_Cao_2022,
Tang_2014,
ThiNgan_ManhToi_2020,
Yang_Feng_deSouza_Wang_2019,
Yang_Li_Lu_2018,
Yuming_Qin_2025_DCDS,
Zang_2022, 
Zhao_Zhu_2015}, and the references therein.  

In the context of the present work, we mention the Velocity-Vorticity formulation of the Navier-Stokes equations, where one considers the velocity equation (in rotational form) and the vorticity equation as independent equations; that is, without enforcing the Biot-Savart coupling $\Omega=\nabla\times U$.  This formulation has become popular in recent years, especially in the regime of vortex-dominated and rotating flows; see, e.g., 
\cite{Benmoussa_Deteix_Yakoubi_2026_JMAA_VVH,G91,GHH90,GHOR15,HOR17,LOR11,LYM06,MF00,OR10,ORS17,WB02,WWW95}.  

In \cite{Larios_Pei_Rebholz_2018}, together with L.~G.~Rebholz, we proposed a Voigt-regularization of the velocity-vorticity formulation of the Navier-Stokes equations (which we called the Velocity-Vorticity-Voigt system, or VVV system), where we put a Voigt regularization term on only the velocity equation, with the idea being to leave the vorticity field in some sense untouched for the sake of preserving the structure of the velocity equation, but recognizing that from the sake of global well-posedness, one still needs some regularization of the velocity equation due to the long-standing difficulties with the 3D Navier-Stokes equations.  Since then, \cite{Duong_2022_VVV_memory} studied the VVV system in the context of damping and memory terms, and the work 
\cite{Yue_Wang_2020_CAMWA_VVV_Attractors_damping} has shown that the VVV system with damping in the vorticity equation has a global attractor, and provided estimates on its dimension. 

In the present work, we extend the VVV system to the case of the 3D MHD system, which we call the VVV-MHD system.  Note that it is not immediately clear how to do this.  For example, the VVV system uses two evolution equations: the velocity and vorticity equation.  Since the MHD system also has a magnetic variable $B$, by considering the current density $J:=\nabla\times B$, it might seem that one would need evolution equations for both of these quantities as well, for a total of four equations.  However, we were pleased to discover that one can prove global well-posedness for the system with only three evolution equations: no additional evolution equation or Voigt regularization for the (analogue of) the current density is needed, and hence the structure of this part of the equations is preserved. This is one of the central results of this paper.  Hence, we propose system \eqref{VVV_MHD} below.  Note that system \eqref{VVV_MHD} is somewhat analogous to the Voigt-MHD (V-MHD) system proposed and studied in \cite{Larios_Titi_2010_MHD}, where the Voigt term is only on the velocity equation.  However, taking the curl of the velocity equation in the V-MHD system yields a vorticity equation that has a Voigt term in it, whereas the such a term is absent in the VVV-MHD equations we study here, so in some sense the VVV-MHD system has lesser regularization than the V-MHD system.  

\begin{remark}
Note that in \cite{Larios_Titi_2010_MHD}, the authors considered the V-MHD equations in the inviscid resistive case ($\nu=0$, $\eta>0$).  We suffer a similar difficulty in both the VVV and VVV-MHD equations, since the vorticity equation has precisely the same structure as the magnetic equation, but here, this means that the vorticity equation seems to require non-zero diffusion.  One could perhaps consider the VVV or the VVV-MHD equations with $\nu=0$ in the velocity equation and $\nu>0$ in the vorticity equation, but this means that the mismatch between the equations would occur not just at the regularization level (due to the Voigt term), but also at the level of the physical model.  In short: if there is to be viscosity in the vorticity equation, there should also be viscosity in the velocity equation.  A similar remark was made in \cite{Larios_Pei_Rebholz_2018}.  Hence, we consider case $\nu>0$ in the velocity and vorticity equations, and as in \cite{Larios_Titi_2010_MHD}, we also require $\eta>0$ in the magnetic equation.
\end{remark}

In addition to proving global well-posedness of the VVV-MHD system, we also prove that as $\alpha\rightarrow0$, the solutions converge to solutions of the the original MHD system in the norm $L^\infty(0,T;L^2)$, where $T>0$ is such that strong solutions to the original MHD system exist and are unique on $[0,T]$.  As usual in other works on the Voigt regularization (see, e.g. \cite{Khouider_Titi_2008,
Larios_Lunasin_Titi_2015,
Larios_Pei_Rebholz_2018,
Pei_2021}), this convergence together with an appropriate ``$\alpha$-energy'' balance gives rise to a blow-up criterion for the original equation.  This form of blow-up criterion was first discovered in \cite{Khouider_Titi_2008} in the context of the surface quasi-geostrophic (SQG) equation.

The Voigt regularization for the 2D MHD equations was first proposed and studied in \cite{Oskolkov_1975_MHD}, and later for the 3D MHD equations in \cite{Larios_Titi_2009,Larios_Titi_2010_MHD} (see also \cite{Catania_2009,Catania_Secchi_2009}).  Recent works on the Voigt regularization for MHD include \cite{Yang_Huang_He_2025_JSC_Voigt_MHD_Elsasser}.

Consider the domain $\mathbb{T}^3$ the three-dimensional periodic space $\mathbb{R}^3/\mathbb{Z}^3 = [0, 1]^3$. 
For $T>0$, the 3D MHD system with full fluid viscosity $\nu> 0$ and magnetic resistivity $\eta> 0$ over $\nT^3\times [0, T)$, is given by
\begin{equation}
       \left\{
       \begin{aligned}
       &\frac{\partial U}{\partial t} 
       - \nu\Delta U 
       + (U\cdot\nabla) U 
       + \nabla P 
       = 
       (B\cdot\nabla) B,
     \\&
       \frac{\partial B}{\partial t} 
       - \eta\Delta B 
       + (U\cdot\nabla) B
       = 
       (B\cdot\nabla) U, 
    \\&
    \nabla \cdot U = 0 = \nabla \cdot B,
    \end{aligned}
    \right.
    \label{MHD}
  \end{equation}
with divergence-free initial data $U(x, 0)=U_0$ and $B(x, 0)=B_0 $;
and the scalar $P=P(x, t)$ represents the unknown pressure. Using the vector identity 
\begin{align}\label{vec_identity}
v\cdot\nabla v = (\nabla\times v)\times v + \tfrac12\nabla|v|^2
\end{align}
and denoting $J:=\nabla\times B$ and $\Omega:=\nabla\times U$,
we can write the first equation in system \eqref{MHD} in the form
\begin{align}\label{curl_form}
    \frac{\partial U}{\partial t} 
       - \nu\Delta U 
       + \Omega\times U 
       + \nabla \Pi 
       = 
       J\times B,
\end{align}
where $\Pi := p + \frac12|U|^2-\frac12|B|^2$.  Also, taking the curl of the first equation in \eqref{MHD}, we obtain the vorticity equation
\begin{align}
    \frac{\partial \Omega}{\partial t} 
       - \nu\Delta \Omega
       + (U\cdot\nabla) \Omega
       = (\Omega\cdot\nabla) U
       + \nabla\times(J\times B),
\end{align}

With the same parameters, we now give below the 3D MHD system with a reformulated vorticity equation and regularization.

\begin{subequations}\label{VVV_MHD}
\begin{empheq}[left=\empheqlbrace]{alignat=3}
\label{VVV_MHD_u}
&(I - \alpha^2\Delta)\frac{\partial u}{\partial t} 
       - \nu\Delta u 
       + w\times u 
       + \nabla p 
       = 
       (b\cdot\nabla) b,
     \\&
\label{VVV_MHD_w}
\frac{\partial w}{\partial t} 
       - \nu\Delta w 
       + (u\cdot\nabla) w
       - (w\cdot\nabla) u  
       = 
       \nabla\times(j\times b),
     \\&
\label{VVV_MHD_b}
       \frac{\partial b}{\partial t} 
       - \eta\Delta b 
       + (u\cdot\nabla) b
       = 
       (b\cdot\nabla) u, 
    \\&
\label{VVV_MHD_div_free}
    0 = \nabla \cdot u = \nabla \cdot w = \nabla \cdot b,
\end{empheq}
\end{subequations}
  
where constant $\alpha>0$, and $I$ denotes the identical operator, and $j=\nabla\times b$ represents the electric current density. 

\subsection{Discussion of some possible variations on the VVV-MHD system}

We note that it should be possible to consider the 3D \eqref{VVV_MHD} equation modified so that a Voigt term (i.e., $-\alpha_M^2\Delta b_t$) also appears in the magnetic equation (and also the velocity equation, but not the $w$ equation).  This should allow one to prove global well-posedness even in the case $\eta=0$ (but $\nu>0$, due to the $w$ equation). Note that a case similar to this for the 3D Voigt-MHD equation (with $\alpha>0$, $\alpha_M>0$, $\nu=\eta=0$) was originally studied in \cite{Larios_Titi_2009}, but without the Velocity-Vorticity split seen in the present work.  We do not study this case here, but leave it as a subject for future work.

We note that one of the authors proposed and studied a VVV modification of the 3D Boussinesq system of ocean dynamics in \cite{Pei_2021}.  Along these lines, one could also include a similar thermal coupling in \eqref{VVV_MHD}, but for the sake of brevity, we omit this in the present work, and leave it as a subject for future work.

\subsection{Organization of the present work}

This paper is organized as follows. In Section~\ref{main}, we state all the major results of this work, and in Section~\ref{pre_main}, we provide without proof all the key theorems and lemmas in fluid dynamics that are essential to our results. Then, in Section~\ref{well-posedness} and \ref{conv}, we show the full proof of the existence and uniqueness results on System~\eqref{VVV_MHD}, as well as the two convergence theorems concerning both Systems~\eqref{MHD} and \eqref{VVV_MHD}, respectively. 

\section{Main Results}\label{main}

We state in this section the major theorems and results of this paper. 

\begin{definition}\label{def-weak}
    For arbitrary $T>0$, with $u_0\in V$, $w_0, b_0\in H$, we call the triplet $(u, w, b)$ a weak solution to System~\eqref{VVV_MHD} on the time interval $[0, T]$, if $u\in C(0, T; V)$, $u_t\in L^2(0, T; V)$, $w, b\in C_{w}(0, T; H)\cap L^2(0, T; V)$, $w_t, b_t\in L^2(0, T; H^{-1})$. Further, $(u, w, b)$ satisfies System~\eqref{VVV_MHD} in the following weak sense:
    \begin{subequations}\label{VVV_MHD_Weak}
      \begin{empheq}[left=\empheqlbrace]{alignat=3}
         &\alpha^2((u_t, \psi)) 
          + (u_t, \psi)
          + \nu((u, \psi))
          + \langle w\times u, \psi\rangle
          = 
          \langle B(b, b), \psi\rangle, \nonumber
          \\&
          \langle w_t, \psi\rangle  
          + \nu((w, \psi)) 
          - \langle B(u, \psi), w \rangle 
          - \langle B(u, u), \psi \rangle   
          = 
          -\langle j\times b, \nabla\times\psi\rangle , \nonumber
          \\&
          \langle b_t, \psi\rangle  
          + \eta((b, \psi)) 
          - \langle B(u, \psi), b\rangle 
          = 
          \langle B(b, u), \psi\rangle , \nonumber
      \end{empheq}
    \end{subequations}
    for all $\psi\in L^2(0, T; V)$. 
\end{definition}

\begin{definition}\label{def-strong}
    For arbitrary $T>0$, with $u_0, w_0, b_0\in V$, we call the triplet $(u, w, b)$ a strong solution to System~\eqref{VVV_MHD} on the time interval $[0, T]$, if it is a weak solution as in Definition~\ref{def-weak}, and further satisfies $w, b\in C(0, T; V)\cap L^2(0, T; D(A))$ and $w_t, b_t\in L^2(0, T; H)$. 
\end{definition}

We state the following well-known local well-posedness result for the 3D MHD system \eqref{MHD}, which can be found in, e.g., \cite{Duvaut_Lions_1972,Sermange_Temam1983,Temam_1997_IDDS}.

\begin{theorem}[Short-time Existence]\label{short-time}
For initial conditions $U_0, B_0\in H^3\cap V$, there exists some time $T>0$, depending on the initial data and on $\nu$ and $\eta$, such that System~\eqref{MHD} possesses  a unique solution $U(t), B(t)\in C_{w}(0, T; H^3)\cap L^{2}(0, T; H^4)$.  
\end{theorem}

Our first major result is the global existence of a unique regular solution to System~\eqref{VVV_MHD}. 
\begin{theorem}[Weak Solution]\label{thm-L2}
For $u_0\in V$ and $w_0, b_0\in H$, System~\eqref{VVV_MHD} has a unique global weak solution in the sense of Definition~\ref{def-weak} which satisfies $\nabla\cdot w = 0$. Further, the following energy identity holds 
    \begin{align}\label{Energy}
        &
        \alpha^2\normL{\nabla u(t)}^2
        +
        \normL{u(t)}^2
        +
        \normL{b(t)}^2
        +
        2\int_{0}^{t}\big(\nu\normL{\nabla u(s)}^2 + \eta\normL{\nabla b(s)}^2\big)\,ds
        \\& \nonumber
        =
        \alpha^2\normL{\nabla u_0}^2
        +
        \normL{u_0}^2
        +
        \normL{b_0}^2.
    \end{align}
\end{theorem}

\begin{theorem}[Strong Solution]\label{thm-regularity}
For $u_0, w_0, b_0\in V$, System~\eqref{VVV_MHD} has a unique global strong solution in the sense of Definition~\ref{def-strong}. If we further assume $u_0, w_0, b_0\in H^s\cap V$ for $\mathbb{N}\ni s\geq 2$, then, the solution satisfies $u\in C_{w}(0, T; H^s\cap V)$ and $w, b\in C_{w}(0, T; H^{s}\cap V)\cap L^2(0, T; H^{s+1}\cap V)$. 
\end{theorem}

The next two theorems concerns the convergence of the solution to System~\eqref{VVV_MHD} with Voigt-regularization term to that of the original MHD system (\eqref{MHD}), on the time-interval of existence of solutions of the latter; as well as the convergence of reformulated $w$ to the original vorticity $\omega$. 

\begin{theorem}[Convergence of $w$ to $\omega=\nabla\times u$ in (\eqref{VVV_MHD}) ]\label{thm-conv1}
    Denote by $\omega:=\nabla\times u$, the vorticity of the formulated flow $u$ and given initial data $u_0, b_0\in H^3\cap V$ and $w_0=\nabla\times u_0$, we have 
    \begin{align*}
      &
      \|\omega(t) - w(t)\|_{L^2}^2
      +
      \alpha^2\|\nabla\omega(t) - \nabla w(t)\|_{L^2}^2
      +
      \int_{0}^{t}\|\nabla\omega(s) - \nabla w(s)\|_{L^2}^2\,ds
      \\&
      \leq
      K_0\alpha^2\big(e^{Ct} - 1\big),
    \end{align*}
    where the constant $K_0$ and $C$ is to be specified in the proof, and where the time-interval $[0, T]$ is such $t<T$ and on which the solution $u(t)\in L^{\infty}(0, T; H^4)\cap L^{2}(0, T; H^4)$ and $w(t)\in L^{\infty}(0, T; H^2)\cap L^{2}(0, T; H^3)$ exists. 
    Namely, as $\alpha\to 0^{+}$, we have $\|\omega(t) - w(t)\|_{L^2}\to 0^{+}$, for a.e. $0<t<T$. 
\end{theorem}

\begin{theorem}[Convergence of the solution of (\eqref{VVV_MHD}) to that of (\eqref{MHD})]\label{thm-conv2}
    For initial conditions $U_0=u_0\in H^3\cap V$, $B_0=b_0\in H^3\cap V$, and $w_0=\nabla\times u_0$, and by denoting $\Omega=\nabla\times U$, and $J=\nabla\times B$, $j=\nabla\times b$, it holds that
    \begin{align*}
      &
      \|u(T)-U(T)\|_{L^2}^2 + \alpha^2\|\nabla u(T) - \nabla U(T)\|_{L^2}^2 + \|w(T) - \Omega(T)\|_{L^2}^2 
      \\&\quad
      + \|b(T) - B(T)\|_{L^2}^2 + \|j(T) - J(T)\|_{L^2}^2
      \\&\quad
      +
      \nu\int_{0}^{T} \big(\|\nabla u(t) - \nabla U(t)\|_{L^2}^2 + \|\nabla w(t) - \nabla\Omega(t)\|_{L^2}^2\big)\,dt
      \\&\quad\quad
      +
      \eta\int_{0}^{T} \big(\|\nabla b(t) - \nabla B(t)\|_{L^2}^2 + \|\nabla j(t) - \nabla J(t)\|_{L^2}^2\big)\,dt
      \\&
      \leq
      C_{T}\alpha^2,
    \end{align*}
    where on the time-interval $[0, T]$, solutions 
    $$u, U, b, B\in L^{\infty}(0, T; H^3)\cap L^{2}(0, T; H^4)$$ 
    and $w(t)\in L^{\infty}(0, T; H^2)\cap L^{2}(0, T; H^3)$ all exist; and where the constant $C$ depends on the $H^3$-norms of $U$, $B$ but not on $\alpha$. 
    Namely, as $\alpha\to 0^{+}$, we have $\|u(t) - U(t)\|_{L^2}\to 0^{+}$, $\|w(t) - \Omega(t)\|_{L^2}\to 0^{+}$, $\|b(t) - B(t)\|_{L^2}\to 0^{+}$, and $\|j(t) - J(t)\|_{L^2}\to 0^{+}$, for a.e. $0<t<T$.
\end{theorem}

Next, we have the following blow-up criterion.  Such a criterion was first stated and proved in the context of the inviscid surface quasi-geostrophic (SQG) equations in \cite{Khouider_Titi_2008} (see also, e.g., \cite{Larios_Titi_2009,Larios_Titi_2010_MHD,Brzezniak_Larios_Safarik_2025_JMFM} and the references therein).

\begin{theorem}[Blow-up Criterion I]\label{thm-blowup}  Under the same conditions and with the same notations as in Theorem~\ref{thm-conv2}, suppose that there exists some $0<T<\infty$ and $\epsilon>0$ such that it holds 
\begin{equation}\label{blow-up_condition}
    \sup\limits_{t\in[0, T]}\limsup\limits_{\alpha\to 0^{+}}\alpha\normL{\nabla u(t)}
    \geq
    \epsilon.
\end{equation}
Then, the solutions to system \eqref{MHD} develop a singularity on $[0, T]$, in the sense that the strong solution ceases to exist. 
\end{theorem}

The following stronger criterion can be proved by following line-by-line an argument in \cite{Larios_Petersen_Titi_Wingate_2015}.

\begin{corollary}[Blow-up Criterion II]\label{coro-blowup}  Under the same conditions and with the same notations as in Theorem~\ref{thm-conv2}, suppose that there exists some $0<T<\infty$ and $\epsilon>0$ such that it holds 
\begin{equation}\label{blow-up_condition2}
    \limsup\limits_{\alpha\to 0^{+}}\alpha\sup\limits_{t\in[0, T]}\normL{\nabla u(t)}
    \geq
    \epsilon.
\end{equation}
Then, the solutions to system \eqref{MHD} develop a singularity on $[0, T]$, in the sense that the strong solution ceases to exist. 
\end{corollary}

\section{Preliminaries}\label{pre_main}

In this section, we first provide the essential preliminaries that we need; then, we list the main theorems and the organization of this paper.

Throughout this paper, we denote by $f_{x_1}$ the partial derivative of $f(x_1,x_2)$ with respect to $x_1$, and likewise for $g_{x_2}$, etc.. We also denote the Lebesgue and the Sobolev spaces by $L^p$ for $0\leq p \leq \infty$ and $H^s = W^{s,2}$ with $s>0$, respectively. Let $\mathcal V$ be the set of all $L$-periodic trigonometric polynomials form $\nR^3$ to $\nR^3$ that are divergence free with zero average. We denote by $H$, $V$, and $V'$ the closures of $\mathcal V$ in the $L^2(\nT^3)^3$ and $H^1(\nT^3)^3$ norms, and the dual space of $V$, respectively, with inner products on $H$ and $V$ as
$$(u, v) = \sum_{i=1}^3\int_{\mathbb{T}^3}u_{i}v_{i}\,dx \text{ \,\,and\,\,  } (\nabla u, \nabla v) = \sum_{i, j=1}^3\int_{\mathbb{T}^3}\partial_{j}u_{i}\partial_{j}v_{i}\,dx,$$
respectively, associated with the norms 
$$\| u \|_{H}=(u, u)^{1/2} \,\,\text{and}\,\, \| u \|_{V}=(\nabla u, \nabla u)^{1/2},$$ where 
$$u(x) = (u_1(x_1, x_2, x_3), u_2(x_1, x_2, x_3), u_3(x_1, x_2, x_3)),$$
and
$$v(x) = (v_1(x_1, x_2, x_3), v_2(x_1, x_2, x_3), v_3(x_1, x_2, x_3)).$$ 
For the sake of convenience, we use $\normL{u}$ and $\normH{u}$ to denote the above norms in $H$ and $V$, respectively. 

We define the Stokes operator $A:= -P_{\sigma}\Delta$ 
with domain $\mathcal{D}(A):= H^2\cap V$, 
where $P_{\sigma}$ is the orthogonal Leray-Helmholtz projector from $L^2$ to $H$. 
Notice that on our domain $\nT^2$, 
we have $A = -\Delta P_{\sigma}$.
Moreover, the Stokes operator can be extended 
to a linear operator from $V$ to $V'$ as 
$$\left<Au, v\right> = (\nabla u, \nabla v) \text{  for all  } v\in V.$$
Moreover, we have the well-known fact that $A^{-1} : H \hookrightarrow \mathcal{D}(A)$ 
is a positive-definite, self-adjoint, and compact operator from $H$ into itself,
thus, $A^{-1}$ possesses an orthonormal basis of positive eigenfunctions $\{ w_{k}\}_{k=1}^{\infty}$ in $H$, corresponding to a sequence of non-increasing sequence of eigenvalues. 
Therefore, $A$ has non-decreasing eigenvalues $\lambda_{k}$, 
i.e., $0 \leq \lambda_1 \leq \lambda_2, \dots$. 
Then, we have the following Poincar\'e inqualities:
\begin{align}\label{Poincare}
    \lambda_1\normL{u}^2
    \leq
    \normH{u}^2 \text{\,\, for\,\,} u\in V;\quad
    \lambda_1\normH{u}^2
    \leq
    \normL{Au}^2 \text{\,\, for\,\,} u\in D(A).
\end{align}
Thus, $\normL{u}$ is equivalent to $\normH{u}$. 

We frequently use the following Agmon's inequality in three-dimensional space.
\begin{equation}\label{agmon1}
    \infnorm{u}^2 \leq c\normH{u} \Htwonorm{u},
\end{equation}
for all $u\in H^2$, as well as 
\begin{equation}\label{agmon2}
    \infnorm{u}^2 \leq c\normL{u}^{1/2} \Htwonorm{u}^{3/2}.
\end{equation}

Also, we state the following Ladyzhenskaya inequality 
\begin{equation}\label{ladyzhenskaya}
    \fournorm{u}^2 \leq c\normL{u}^{1/2} \normH{u}^{3/2},
\end{equation}
for all $u\in V$, which is a special case of the following interpolation result 
that is frequently used in this paper (see, e.g., \cite{Nirenberg_1959_AnnPisa} for a detailed proof).
Assume $1 \leq q, r \leq \infty$, and $0<\gamma<1$.  
For $v\in L^q(\mathbb{T}^{n})$, such that  $\partial^\alpha v\in L^{r} (\mathbb{T}^{n})$, for $|\alpha|=m$, then 
\begin{align}\label{PT1}
\|\partial_{s}v\|_{L^{p}} \leq C\|\partial^{\alpha}v\|_{L^{r}}^{\gamma}\| v\|_{L^{q}}^{1-\gamma},
\,\,\,\text{where}\,\,\,
\frac{1}{p} - \frac{s}{n} = \left(\frac{1}{r} - \frac{m}{n}\right) \gamma+ \frac{1}{q}(1-\gamma).
\end{align}

The following lemma, regarding the bilinear term 
\begin{align*}
     B(u, v) := P_{\sigma}((u\cdot\nabla)v),\quad (u, v \in \mathcal{V}),
\end{align*}
which can be extended to a continuous map (with a bit abuse of notation) 
$B : V \times V \to V'$,
\begin{align*}
     \left<B(u, v), w\right> = \int_{\mathbb{T}^2}(u\cdot\nabla v)\cdot w\,dx, \quad (u, v, w\in \mathcal{V}),
\end{align*}
is frequently used in this article (See, e.g., \cite{Constantin_Foias_1988,Temam_2001_Th_Num,Foias_Manley_Rosa_Temam_2001} for details).
\begin{lemma}
    \begin{subequations}
      \begin{align}
      \label{symm1}
      \ip{B(u,v)}{w}_{V'} &= -\ip{B(u,w)}{v}_{V'}, 
      \quad\forall\;u, v, w\in V,\\
      \label{symm2}
      \ip{B(u,v)}{v}_{V'} &= 0,
      \quad\forall\;u, v, w\in V.
    \end{align}
  \end{subequations}
  \begin{subequations}
    \begin{align}
      \label{B:424}
      |\ip{B(u,v)}{w}_{V'}|
      &\leq 
      C\normL{u}^{1/2} \normH{u}^{1/2} \normH{v} \normL{w}^{1/2} \normH{w}^{1/2},
      \\
      \label{B:442}
      |\ip{B(u,v)}{w}_{V'}|
      &\leq 
      C\normL{u}^{1/2} \normH{u}^{1/2} \normH{u}^{1/2} \normL{Av}^{1/2} \normL{w},
      \\
      \label{B:inf22}
      |\ip{B(u,v)}{w}_{V'}|
      &\leq 
      C\normL{u}^{1/2} \normL{Au}^{1/2} \normH{v} \normL{w},
    \end{align}
for all $u$, $v$, $w$ in the largest spaces $H$, $V$, or $D(A)$, for which the right-hand sides of the inequalities are finite.
  \end{subequations}
  \begin{align}\label{enstrophy_miracle}
    \ip{B(w, w)}{Aw} = 0,\; \forall w\in D(A), \quad \text{(in two-dimension)}
\end{align}
  and the Jacobi identity   
  \begin{align}\label{jacobi}
    \ip{B(u, w)}{Aw} + \ip{B(w, u)}{Aw} +\ip{B(w, w)}{Au} = 0.
  \end{align}
\end{lemma}

For the sake of completeness, we state the following uniform Gr\"{o}nwall's inequality, proved in \cite{Jones_Titi_1992} (see also \cite{Farhat_Lunasin_Titi_2017_Horizontal} and the references therein), which will be used frequently throughout the paper. 
\begin{lemma}
\label{Gronwall}
Suppose that $Y(t)$ is a locally integrable and absolutely continuous function that satisfies the following: 
$$\frac{d Y}{d t} + \alpha(t) Y \leq \beta(t), \quad\text{ a.e. on } (0, \infty), $$
such that 
$$\liminf_{t \to \infty} \int_{t}^{t+\tau} \alpha(s)\,ds \geq \gamma, \quad\quad\quad \limsup_{t \to \infty} \int_{t}^{t+\tau} \alpha^{-}(s)\,ds < \infty, $$
and 
$$\lim_{t \to \infty} \int_{t}^{t+\tau} \beta^{+}(s)\,ds = 0, $$
for fixed $\tau > 0$, and $\gamma > 0$, 
where 
$\alpha^{-} = \max\{-\alpha, 0\}$ 
and 
$\beta^{+} = \max\{\beta, 0\}$. 
Then, $Y(t) \to 0$ at an exponential rate as $t \to \infty$.
\end{lemma}

\section{Proof of Well-posedness Results}\label{well-posedness}
For the sake of simplicity, we provide only the {\it a priori} estimates in order to obtain the global well-posedness of the solution to system \eqref{VVV_MHD}, following the approach in \cite{Larios_Pei_2017_MHDB_PS} and the reference therein at the formal level. The argument can be made rigorous by the Galerkin approximation method, of which the details we omit here. We start with the $H^1$ estimates of $u$ and $L^2$-estimates of $b$. 

\subsection{\texorpdfstring{$H^1$-estimates of $u$ and $L^2$-estimates of $b$ in Theorem~\ref{thm-L2}}{}}\label{proof_L2}
Taking the inner product of equations \eqref{VVV_MHD_u} and \eqref{VVV_MHD_b} in System~\eqref{VVV_MHD} with $u$ and $b$, respectively, we obtain 
\begin{align*}
    &
    \frac{1}{2}\frac{d}{dt} \big( \norm{u}_{L^2}^2 
    + 
    \alpha^2\norm{\nabla u}_{L^2}^2 
    +
    \norm{b}_{L^2}^2 \big)
    + 
    \nu \norm{\nabla u}_{L^2}^2 
    +
    \eta \norm{\nabla b}_{L^2}^2
    \\&
    = 
    \int_{\mathbb{T}^3}(b \cdot \nabla)b\cdot u
    + 
    \int_{\mathbb{T}^3}(b \cdot \nabla)u\cdot b
    =0,
\end{align*}
where we integrated by part on the right side of the above equation and used the divergence-free condition on $u$ and $b$. 
Then, integrating in time on $[0, t]$ where $t\in [0, T]$ yields the desired energy identity \eqref{Energy}. 

\subsection{\texorpdfstring{$H^1$}--estimates of Theorem~\ref{thm-regularity}}\label{b-H1}
In this section, we provide $H^1$ estimates of $b$. 
Taking the curl of \eqref{VVV_MHD_b}, we obtain
\begin{align}\label{eq_curl_b}
    \frac{\partial j}{\partial t} 
    - 
    \eta \triangle j 
    = 
    \nabla \times \big((b \cdot \nabla) u - (u\cdot \nabla)b\big)
\end{align}
where $j := \nabla \times b$. 
Multiplying the above equation by $j$ and integrating by parts, we have 
\begin{align}\label{j_energy}
    \frac{1}{2}\frac{d}{dt} \norm{j}_{L^2}^2 
    + 
    \eta \norm{\nabla j}_{L^2}^2 
    = 
    -\int_{\mathbb{T}^3}(b\cdot \nabla)u\cdot(\nabla \times j) 
    + 
    \int_{\mathbb{T}^3}(u\cdot \nabla)b\cdot(\nabla \times j). 
\end{align}
The first integral on the right side of the above equations is bounded by 
\begin{align*}
    \left|\int_{\mathbb{T}^3}(b\cdot \nabla)u\cdot(\nabla \times j) \right| 
    & 
    \leq 
    \norm{b}_{L^\infty}\norm{\nabla u}_{L^2}\norm{\nabla \times j}_{L^2}
    \leq 
    C \norm{j}_{L^2}^{1/2}\norm{\nabla j}_{L^2}^{3/2}\norm{\nabla u}_{L^2}^2
    \\& 
    \leq 
    C_\alpha \norm{j}_{L^2}^{1/2}\norm{\nabla j}_{L^2}^{3/2}
    \leq 
    C_{\alpha, \eta}\norm{j}_{L^2}^2 + \frac{3\eta}{4}\norm{\nabla j}_{L^2}^2
\end{align*}
where we applied Agmon's inequality \eqref{agmon1} and Young's inequality. 
As for the second integral, we have 
\begin{align*}
    \left|\int_{\mathbb{T}^3}(u\cdot \nabla)b\cdot(\nabla \times j)\right|
    & 
    \leq 
    \norm{u}_{L^6}\norm{\nabla b}_{L^3}\norm{\nabla \times j}_{L^2}
    \leq 
    C\norm{\nabla u}_{L^2}\norm{j}_{L^2}^{1/2}\norm{\nabla j}_{L^2}^{3/2}
    \\&
    \leq 
    C_{\alpha, \eta}\norm{j}_{L^2}^2 + \frac{3\eta}{4}\norm{\nabla j}_{L^2}^2
\end{align*}
Hence, we obtain 
\begin{align}\label{b_H1_est}
    \frac{1}{2}\frac{d}{dt}\norm{j}_{L^2}^2 
    + 
    \frac{\eta}{4}\norm{\nabla j}_{L^2}^2 
    \leq 
    C_{\alpha, \eta}\norm{j}_{L^2}^2. 
\end{align}
Using Gr\"onwall's inequality, we obtain 
\begin{align}\label{h1_h2_b}
    b \in L^\infty((0,\infty); H^1)\cap L^2((0,T); H^2).\quad\text{\qedsymbol}
\end{align}

\subsection{\texorpdfstring{$L^2$}--estimates of \texorpdfstring{$w$}-}\label{w-L2}
Now, we are ready to provide the $L^2$-estimates of $w$. 
Taking the inner-product of \eqref{VVV_MHD_w} with $w$, integrating by parts, and adding the resulted equation to \eqref{j_energy}, we have 
\begin{align*}
    &
    \frac{1}{2} \frac{d}{dt} \big(\norm{w}_{L^2}^2 + \norm{j}_{L^2}^2\big)
    + 
    \nu \norm{\nabla w}_{L^2}^2 
    +
    \eta \norm{\nabla j}_{L^2}^2 
    \\&
    = 
    -\int_{\mathbb{T}^3}(b\cdot \nabla)u\cdot(\nabla \times j) 
    + 
    \int_{\mathbb{T}^3}(u\cdot \nabla)b\cdot(\nabla \times j)
    \\&\quad
    +
    \int_{\mathbb{T}^3}(w \cdot \nabla) u\cdot w 
    + 
    \int_{\mathbb{T}^3}(\nabla \times (j\times b))\cdot w.
\end{align*}
We have already estimated the first two integrals on the right side of the above equation. In order to bound the third integral, we apply Sobolev, Agmon's, and Young's inequalities and get 
\begin{align*}
    \left|\int_{\mathbb{T}^3}(w \cdot \nabla) u\cdot w\right| 
    & 
    \leq 
    \norm{w}_{L^3} \norm{\nabla u}_{L^2} \norm{w}_{L^6} 
    \leq
    C \norm{w}_{L^2}^{1/2} \norm{\nabla u}_{L^2}\norm{\nabla w}_{L^2}^{3/2} 
    \\& 
    \leq 
    C_\nu \norm{\nabla u}_{L^2}^4 \norm{w}_{L^2}^2
    +
    \frac{\nu}{4} \norm{\nabla w}_{L^2}^{2} 
    \\&
    \leq
    C_{\alpha, \nu} \norm{w}_{L^2}^2
    +
    \frac{\nu}{4} \norm{\nabla w}_{L^2}^{2}. 
\end{align*}
As for the fourth integral, we have:
\begin{align*}
    \left|\int_{\mathbb{T}^3}(\nabla \times (j\times b)\cdot w\right|
    & 
    = 
    \left|\int_{\mathbb{T}^3}(j\times b)\cdot (\nabla\times w)\right| 
    \leq 
    \norm{b}_{L^6} \norm{j}_{L^3} \norm{\nabla \times w}_{L^2} 
    \\& 
    \leq 
    C\norm{j}_{L^2}^{3/2} \norm{\nabla j}_{L^2}^{1/2} \norm{\nabla w}_{L^2} 
    \leq 
    C_{\alpha, \eta}\norm{j}_{L^2}^{1/2} \norm{\nabla j}_{L^2}^{1/2} \norm{\nabla w}_{L^2} 
    \\&
    \leq
    C_{\alpha,\nu,\eta}\norm{j}_{L^2} \norm{\nabla j}_{L^2}
    +
    \frac{\nu}{4}\norm{\nabla w}_{L^2}^2
    \\&
    \leq 
    C_{\alpha,\nu,\eta}\norm{j}_{L^2} 
    + 
    \frac{\eta}{4} \norm{\nabla j}_{L^2}^2 
    +
    \frac{\nu}{4}\norm{\nabla w}_{L^2}^2.
\end{align*}
Combining the above estimates, we obtain
\begin{align}\label{w_L2_est}
    &
    \frac{d}{dt}(\norm{w}_{L^2}^2 + \norm{j}_{L^2}^2)
    +
    \nu\normL{\nabla w}^2
    +
    \eta\normL{\nabla j}^2
    \leq 
    C_{\alpha,\nu,\eta} (\norm{w}_{L^2}^2 + \norm{j}_{L^2}^2).
\end{align}
Now, Gr\"onwall's inequality implies  
\begin{align*}
    w \in L^\infty((0,T); H) \cap L^2(0,T; V).\quad\text{\qedsymbol}
\end{align*}

\subsection{\texorpdfstring{$H^2$-estimates of $u$ and $b$}{}}\label{u-b-higher-reg}
In this section, we provide $H^2$-estimates of $u$ and $b$.  We start with the bound on $\normL{\Delta u}^2$. We first multiply \eqref{VVV_MHD_u} by $-\Delta u$, and \eqref{VVV_MHD_b} by $-\Delta b$, add, and integrate by parts, to obtain
\begin{align*}
    &
    \frac{1}{2} \frac{d}{dt}\big(\norm{\nabla u}_{L^2}^2 + \alpha^2\norm{\Delta u}_{L^2}^2 +  \norm{\nabla b}_{L^2}^2\big)
    + 
    \nu \norm{\Delta u}_{L^2}^2 
    +
    \eta \norm{\Delta b}_{L^2}^2 
    \\&
    = 
    \int_{\mathbb{T}^3}(w\times u)\cdot(\Delta u) 
    - 
    \int_{\mathbb{T}^3}(b\cdot \nabla)b\cdot(\Delta u)
    \\&\quad
    -
    \int_{\mathbb{T}^3}(b \cdot \nabla) u\cdot (\Delta b) 
    + 
    \int_{\mathbb{T}^3}(u \cdot \nabla) b\cdot (\Delta b).
\end{align*}
Then, we apply Sobolev, H\"older's, Agmon's, and Young's inequalities to estimate each of the four integrals on the right side of the above equations. The first one is bounded by 
\begin{align*}
    \int_{\mathbb{T}^3}|w||u||\Delta u|
    &
    \leq
    C\normLp{3}{w}\normLp{6}{u}\normL{\Delta u}
    \\&
    \leq
    C\normL{w}^{1/2}\normL{\nabla w}^{1/2}\normL{\nabla u}\normL{\Delta u}
    \\&
    \leq
    C\normL{\nabla w}\normL{\nabla u}\normL{\Delta u}
    \\&
    \leq
    C_{\nu}\normL{\nabla w}^2\normL{\nabla u}^2
    +
    \frac{\nu}{8}\normL{\Delta u}.
\end{align*} 
The second integral is bounded by 
\begin{align*}
    \int_{\mathbb{T}^3}|b||\nabla b||\Delta u|
    &
    \leq
    C\normLp{6}{b}\normLp{3}{\nabla b}\normL{\Delta u}
    \leq
    C\normL{\nabla b}^{3/2}\normL{\Delta b}^{1/2}\normL{\Delta u}
    \\&
    \leq
    C_{\nu}\normL{\nabla b}^3\normL{\Delta b}
    +
    \frac{\nu}{8}\normL{\Delta u}^2
    \\&
    \leq
    \tilde{C}_{\nu}\normL{\nabla b}\normL{\Delta b}
    +
    \frac{\nu}{8}\normL{\Delta u}^2
    \\&
    \leq
    \tilde{C}_{\nu}\normL{\nabla b}^2
    +
    \frac{\eta}{8}\normL{\Delta b}^2
    +
    \frac{\nu}{8}\normL{\Delta u}^2.
\end{align*}
Next, we bound the third integral by
\begin{align*}
    \int_{\mathbb{T}^3}|b||\nabla u||\Delta b|
    &
    \leq
    C\normLp{6}{b}\normLp{3}{\nabla u}\normL{\Delta b}
    \\&
    \leq
    C\normL{\nabla b}\normL{\nabla u}^{1/2}\normL{\Delta u}^{1/2}\normL{\Delta b}
    \\&
    \leq
    C_{\eta}\normL{\nabla b}^2\normL{\nabla u}\normL{\Delta u}
    +
    \frac{\eta}{8}\normL{\Delta b}^2
    \\&
    \leq
    C_{\nu,\eta}\normL{\nabla b}^4\normL{\nabla u}^2
    +
    \frac{\nu}{8}\normL{\Delta u}^2
    +
    \frac{\nu}{8}\normL{\Delta b}^2
    \\&
    \leq
    \tilde{C}_{\nu,\eta}\normL{\nabla u}^2
    +
    \frac{\eta}{8}\normL{\Delta u}^2
    +
    \frac{\nu}{8}\normL{\Delta b}^2.
\end{align*}
As for the bound on the last integral, we have 
\begin{align*}
    \int_{\mathbb{T}^3}|u||\nabla b||\Delta b|
    &
    \leq
    C\normLp{6}{u}\normLp{3}{\nabla b}\normL{\Delta b}
    \\&
    \leq
    C\normL{\nabla u}\normL{\nabla b}^{1/2}\normL{\Delta b}^{3/2}
    \\&
    \leq
    C_{\eta}\normL{\nabla u}^4\normL{\nabla b}^2
    +
    \frac{\eta}{8}\normL{\Delta b}^2
    \\&
    \leq
    \tilde{C}_{\eta}\normL{\nabla b}^2
    +
    \frac{\eta}{8}\normL{\Delta u}^2.
\end{align*}
Thus, combining all the above estimates, we obtain
\begin{align*}
   &\frac{d}{dt}\big(\norm{\nabla u}_{L^2}^2 + \alpha^2\norm{\Delta u}_{L^2}^2 +  \norm{\nabla b}_{L^2}^2\big)
   + 
   \nu \norm{\Delta u}_{L^2}^2 
   +
   \eta \norm{\Delta b}_{L^2}^2 
   \\&
   \leq
   C\big(\norm{\nabla u}_{L^2}^2 + \alpha^2\norm{\Delta u}_{L^2}^2 +  \norm{\nabla b}_{L^2}^2\big),
\end{align*}
where the constant $C$ depends on $\nu, \eta$ and $\normL{\nabla w}^2$, which is integrable in time. Hence, Gr\"onwall's inequality completes the proof, i.e., we have $$u\in C(0, T; D(A)).$$

We are ready to estimate the $H^2$-bounds of $b$. To do so, we multiply \eqref{VVV_MHD_u} and \eqref{VVV_MHD_b} by $-\Delta u$ and $\Delta^2 b$, respectively, add, integrate by parts, and obtain
\begin{align*}
    &
    \frac{1}{2} \frac{d}{dt}\big(\norm{\nabla u}_{L^2}^2 + \alpha^2\norm{\Delta u}_{L^2}^2 +  \norm{\Delta b}_{L^2}^2\big)
    + 
    \nu \norm{\Delta u}_{L^2}^2 
    +
    \eta \norm{\nabla\Delta b}_{L^2}^2 
    \\&
    = 
    \int_{\mathbb{T}^3}(w\times u)\cdot(\Delta u) 
    - 
    \int_{\mathbb{T}^3}(b\cdot \nabla)b\cdot(\Delta u)
    \\&\quad
    +
    \int_{\mathbb{T}^3}(b \cdot \nabla) u\cdot (\Delta^2 b) 
    - 
    \int_{\mathbb{T}^3}(u \cdot \nabla) b\cdot (\Delta^2 b).
\end{align*}
Then, we bound each of the four integrals on the right side of the above equation. The first one is estimated in the same way as in the $H^2$-estimates of $u$, so we proceed to bound the second integral by
\begin{align*}
    \int_{\mathbb{T}^3}|b||\nabla b||\Delta u|
    &
    \leq
    C\normLp{3}{b}\normLp{6}{\nabla b}\normL{\Delta u}
    \\&
    \leq
    C\normL{b}^{1/2}\normL{\nabla b}^{1/2}\normL{\Delta b}\normL{\Delta u}
    \\&
    \leq
    C_{\nu}\normL{\Delta b}^2
    +
    \frac{\nu}{8}\normL{\Delta u}^2.
\end{align*}
As for the third integral, we first integrate by parts, then bound it by  
\begin{align*}
    &
    \int_{\mathbb{T}^3}|\nabla b||\nabla u||\nabla \Delta b|
    +
    \int_{\mathbb{T}^3}|b||\Delta u||\nabla \Delta b|
    \\&
    \leq
    C\normLp{3}{\nabla b}\normLp{6}{\nabla u}\normL{\nabla \Delta b}
    +
    C\normLp{\infty}{b}\normL{\Delta u}\normL{\nabla \Delta b}
    \\&
    \leq
    C\normL{\nabla b}^{1/2}\normL{\Delta b}^{1/2}\normL{\Delta u}\normL{\nabla \Delta b}
    \\&
    \leq
    C\normL{\Delta u}\normL{\Delta b}\normL{\nabla \Delta b}
    \\&
    \leq
    C_{\eta}\normL{\Delta u}^2\normL{\Delta b}^2
    +
    \frac{\eta}{4}\normL{\nabla \Delta b}^2
    \\&
    \leq
    \tilde{C}_{\eta}\normL{\Delta b}^2
    +
    \frac{\eta}{4}\normL{\nabla \Delta b}^2,
\end{align*}
where we used the $H^2$-bound of $u$ in the last inequality. As for the last integral, we proceed similarly, i.e., integrate by parts, and bound it by 
\begin{align*}
    &
    \int_{\mathbb{T}^3}|\nabla u||\nabla b||\nabla \Delta b|
    +
    \int_{\mathbb{T}^3}|u||\Delta b||\nabla \Delta b|
    \\&
    \leq
    C\normLp{3}{\nabla u}\normLp{6}{\nabla b}\normL{\nabla \Delta b}
    +
    C\normLp{\infty}{u}\normL{\Delta b}\normL{\nabla \Delta b}
    \\&
    \leq
    C\normL{\nabla u}^{1/2}\normL{\Delta u}^{1/2}\normL{\Delta b}\normL{\nabla \Delta b}
    \\&
    \leq
    C_{\eta}\normL{\nabla u}\normL{\Delta u}\normL{\Delta b}^2
    +
    \frac{\eta}{4}\normL{\nabla \Delta b}^2
    \\&
    \leq
    \tilde{C}_{\eta}\normL{\Delta b}^2
    +
    \frac{\eta}{4}\normL{\nabla \Delta b}^2. 
\end{align*}
Therefore, combining all the above estimates, we obtain
\begin{align*}
   &\frac{d}{dt}\big(\norm{\nabla u}_{L^2}^2 + \alpha^2\norm{\Delta u}_{L^2}^2 +  \norm{\Delta b}_{L^2}^2\big)
   + 
   \nu \norm{\Delta u}_{L^2}^2 
   +
   \eta \norm{\nabla\Delta b}_{L^2}^2 
   \\&
   \leq
   \tilde{C}\big(\norm{\nabla u}_{L^2}^2 + \alpha^2\norm{\Delta u}_{L^2}^2 +  \norm{\Delta b}_{L^2}^2\big),
\end{align*}
where the above constant $\tilde{C}$ depends on $\nu, \eta$, $H^2$-bound of $u$, and $\normL{\nabla w}^2$, which is integrable in time. Hence, Gr\"onwall's inequality completes the proof, i.e., we have $$b\in C(0, T; D(A))\cap L^2(0, T; H^3\cap V).\quad\text{\qedsymbol}$$

\subsection{\texorpdfstring{$H^1$}--estimates of \texorpdfstring{$w$}-}\label{w-H1}
In this section, with $H^2$-bounds on $u$ and $b$ in hand, we are ready to provide $H^1$-estimates of $w$. Multiplying equations \eqref{VVV_MHD_u}, \eqref{VVV_MHD_w}, and \eqref{VVV_MHD_b} by $-\Delta u$ and $-\Delta w$, and $-\Delta b$, respectively, integrating by parts, and adding the three resulted equations together, we have 
\begin{align*}
    &
    \frac{1}{2} \frac{d}{dt}\big(\norm{\nabla u}_{L^2}^2 + \alpha^2\norm{\Delta u}_{L^2}^2 + \norm{\nabla w}_{L^2}^2 + \norm{\nabla b}_{L^2}^2\big)
    \\&\quad
    + 
    \nu \norm{\Delta u}_{L^2}^2 
    +
    \nu \norm{\Delta w}_{L^2}^2 
    +
    \eta \norm{\Delta b}_{L^2}^2 
    \\&
    = 
    \int_{\mathbb{T}^3}(w\times u)\cdot(\Delta u) 
    - 
    \int_{\mathbb{T}^3}(b\cdot \nabla)b\cdot(\Delta u)
    \\&\quad
    +
    \int_{\mathbb{T}^3}(u\cdot\nabla)w\cdot(\Delta w)
    -
    \int_{\mathbb{T}^3}(w\cdot\nabla)u\cdot(\Delta w)
    -
    \int_{\mathbb{T}^3}\nabla\times(j\times b)\cdot(\Delta w)
    \\&\quad\quad
    -
    \int_{\mathbb{T}^3}(b \cdot \nabla) u\cdot (\Delta b) 
    + 
    \int_{\mathbb{T}^3}(u \cdot \nabla) b\cdot (\Delta b).
\end{align*}
The second and the last two integrals on the right side of the above equation are estimated analogously as in the $H^2$-estimates of $u$, so we proceed and bound the first integral by 
\begin{align*}
    \int_{\mathbb{T}^3}|w||u||\Delta u|
    &
    \leq
    C\normLp{\infty}{w}\normL{u}\normL{\Delta u}
    \\&
    \leq
    C\normL{u}\normL{\nabla w}^{1/2}\normL{\Delta w}^{1/2}\normL{\Delta u}
    \\&
    \leq
    C_{\nu}\normL{u}^2\normL{\nabla w}\normL{\Delta w}
    +
    \frac{\nu}{8}\normL{\Delta u}^2
    \\&
    \leq
    C_{\nu}\normL{u}^4\normL{\nabla w}^2
    +
    \frac{\nu}{8}\normL{\Delta w}^2
    +
    \frac{\nu}{8}\normL{\Delta u}^2
    \\&
    \leq
    \tilde{C}_{\nu}\normL{\nabla w}^2
    +
    \frac{\nu}{8}\normL{\Delta w}^2
    +
    \frac{\nu}{8}\normL{\Delta u}^2.
\end{align*} 
For the third integral, we bound it by 
\begin{align*}
    \int_{\mathbb{T}^3}|u||\nabla w||\Delta w|
    &
    \leq
    C\normLp{6}{u}\normLp{3}{\nabla w}\normL{\Delta w}
    \\&
    \leq
    C\normL{\nabla u}\normL{\nabla w}^{1/2}\normL{\Delta w}^{3/2}
    \\&
    \leq
    C_{\nu}\normL{\nabla u}^4\normL{\nabla w}^{2}
    +
    \frac{\nu}{8}\normL{\Delta w}^{2}
    \\&
    \leq
    \tilde{C}_{\nu}\normL{\nabla w}^{2}
    +
    \frac{\nu}{8}\normL{\Delta w}^{2},  
\end{align*}
while the fourth integral is bounded similarly by 
\begin{align*}
    \int_{\mathbb{T}^3}|w||\nabla u||\Delta w|
    &
    \leq
    C\normLp{\infty}{w}\normL{\nabla u}\normL{\Delta w}
    \\&
    \leq
    C\normL{\nabla u}\normL{\nabla w}^{1/2}\normL{\Delta w}^{3/2}
    \\&
    \leq
    \tilde{C}_{\nu}\normL{\nabla w}^{2}
    +
    \frac{\nu}{8}\normL{\Delta w}^{2}.  
\end{align*}
Lastly, we bound the fifth integral by 
\begin{align*}
    &
    \int_{\mathbb{T}^3}|\nabla j||b||\Delta w|
    +
    \int_{\mathbb{T}^3}|j||\nabla b||\Delta w|
    \\&
    \leq
    C\normL{\nabla j}\normLp{\infty}{b}\normL{\Delta w}
    +
    C\normLp{6}{j}\normLp{3}{\nabla b}\normL{\Delta w}
    \\&
    \leq
    C\normL{\nabla b}^{1/2}\normL{\Delta b}^{3/2}\normL{\Delta w}
    \leq
    C_{\nu}\normL{\nabla b}\normL{\Delta b}^{3}
    +
    \frac{\nu}{8}\normL{\Delta w}^2
    \\&
    \leq
    \tilde{C}_{\nu}\normL{\nabla b}\normL{\Delta b}
    +
    \frac{\nu}{8}\normL{\Delta w}^2
    \leq
    \tilde{C}_{\nu,\eta}\normL{\nabla b}^2 
    +
    \frac{\eta}{8}\normL{\Delta b}^2
    +
    \frac{\nu}{8}\normL{\Delta w}^2, 
\end{align*}
where we used the $H^2$-bound on $b$ in the last but second inequality above. 
Combining all the above estimates, we obtain 
\begin{align*}
   &\frac{d}{dt}\big(\norm{\nabla u}_{L^2}^2 + \alpha^2\norm{\Delta u}_{L^2}^2 + \norm{\nabla w}^2 + \norm{\nabla b}_{L^2}^2\big)
   \\&\quad
   + 
   \nu \norm{\Delta u}_{L^2}^2 
   +
   \nu \norm{\Delta w}_{L^2}^2
   +
   \eta \norm{\Delta b}_{L^2}^2 
   \\&
   \leq
   C\big(\norm{\nabla u}_{L^2}^2 + \alpha^2\norm{\Delta u}_{L^2}^2 + \norm{\nabla w}^2 + \norm{\nabla b}_{L^2}^2\big),
\end{align*}
where the above constant $C$ depends on $\nu, \eta$, and the $H^2$-bounds on $u$ and $b$. Hence, Gr\"onwall's inequality leads to the desired conclusion, i.e., we have proved $$w\in C(0, T; V)\cap L^2(0, T; H^2\cap V).\quad\text{\qedsymbol}$$

\subsection{\texorpdfstring{$H^3$-estimates of $u$ and $b$}{}}\label{u-b-higher-reg-H3}
In this section, we provide $H^3$-estimates of $u$ and $b$. First, to obtain the $H^3$-bound on $u$, we multiply \eqref{VVV_MHD_u} by $\Delta^2 u$, and \eqref{VVV_MHD_b} by $\Delta^2 b$, respectively, integrate by parts, and add, and get
\begin{align*}
    &
    \frac{1}{2} \frac{d}{dt}\big(\norm{\Delta u}_{L^2}^2 + \alpha^2\norm{\nabla\Delta u}_{L^2}^2 +  \norm{\Delta b}_{L^2}^2\big)
    + 
    \nu \norm{\nabla\Delta u}_{L^2}^2 
    +
    \eta \norm{\nabla\Delta b}_{L^2}^2 
    \\&
    = 
    -\int_{\mathbb{T}^3}(w\times u)\cdot(\Delta^2 u) 
    + 
    \int_{\mathbb{T}^3}(b\cdot \nabla)b\cdot(\Delta^2 u)
    \\&\quad
    +
    \int_{\mathbb{T}^3}(b \cdot \nabla) u\cdot (\Delta^2 b) 
    - 
    \int_{\mathbb{T}^3}(u \cdot \nabla) b\cdot (\Delta^2 b).
\end{align*}
Then, we estimate each of the four integrals on the right side of the above equation. After integration by parts, and applying Sobolev, H\"older's, and Agmon's inequalities, we bound the first integral by
\begin{align*}
    &
    \int_{\mathbb{T}^3}|\nabla w||u||\nabla\Delta u|
    +
    \int_{\mathbb{T}^3}|w||\nabla u||\nabla\Delta u|
    \\&
    \leq
    C\normL{\nabla w}\normLp{\infty}{u}\normL{\nabla\Delta u}
    +
    C\normLp{6}{w}\normLp{3}{\nabla u}\normL{\nabla\Delta u}
    \\&
    \leq
    C\normL{\nabla w}\normL{\nabla u}^{1/2}\normL{\Delta u}^{1/2}\normL{\nabla\Delta u}
    \leq
    C\normL{\Delta u}\normL{\nabla\Delta u}
    \\&
    \leq
    C\normL{\Delta u}^2
    +
    \frac{\nu}{4}\normL{\nabla\Delta u}^2,
\end{align*}
where we also used Poincar\'e's and Young's inequalities in the last two steps. 
Next, the second integral is bounded by 
\begin{align*}
    &
    \int_{\mathbb{T}^3}|\nabla b||\nabla b||\nabla\Delta u|
    +
    \int_{\mathbb{T}^3}|b||\Delta b||\nabla\Delta u|
    \\&
    \leq
    C\normLp{3}{\nabla b}\normLp{6}{\nabla b}\normL{\nabla\Delta u}
    +
    C\normLp{\infty}{b}\normL{\Delta b}\normL{\nabla\Delta u}
    \\&
    \leq
    C\normL{\nabla b}^{1/2}\normL{\Delta b}^{3/2}\normL{\nabla\Delta u}
    \leq
    C\normL{\Delta b}\normL{\nabla\Delta u}
    \\&
    \leq
    C\normL{\Delta b}^2
    +
    \frac{\nu}{4}\normL{\nabla\Delta u}^2,
\end{align*}
where we used the previously-obtained bounds on $\normL{\nabla b}$ and $\normL{\Delta b}$. 
As for the third integral, we proceed similarly and it is bounded by 
\begin{align*}
    &
    \int_{\mathbb{T}^3}|\nabla b||\nabla u||\nabla\Delta b|
    +
    \int_{\mathbb{T}^3}|b||\Delta u||\nabla\Delta b|
    \\&
    \leq
    C\normLp{3}{\nabla b}\normLp{6}{\nabla u}\normL{\nabla\Delta b}
    +
    C\normLp{\infty}{b}\normL{\Delta u}\normL{\nabla\Delta b}
    \\&
    \leq
    C\normL{\nabla b}^{1/2}\normL{\Delta b}^{1/2}\normL{\Delta u}\normL{\nabla\Delta b}
    \leq
    C\normL{\Delta u}\normL{\nabla\Delta b}
    \\&
    \leq
    C\normL{\Delta u}^2
    +
    \frac{\eta}{4}\normL{\nabla\Delta b}^2. 
\end{align*}
The last integral is bounded by 
\begin{align*}
    &
    \int_{\mathbb{T}^3}|\nabla u||\nabla b||\nabla\Delta b|
    +
    \int_{\mathbb{T}^3}|u||\Delta b||\nabla\Delta b|
    \\&
    \leq
    C\normLp{3}{\nabla u}\normLp{6}{\nabla b}\normL{\nabla\Delta b}
    +
    C\normLp{\infty}{u}\normL{\Delta b}\normL{\nabla\Delta b}
    \\&
    \leq
    C\normL{\nabla u}^{1/2}\normL{\Delta u}^{1/2}\normL{\Delta b}\normL{\nabla\Delta b}
    \leq
    C\normL{\Delta b}\normL{\nabla\Delta b}
    \\&
    \leq
    C\normL{\Delta b}^2
    +
    \frac{\eta}{4}\normL{\nabla\Delta b}^2. 
\end{align*}
Thus, combining all the above estimates, we obtain
\begin{align*}
   &\frac{d}{dt}\big(\norm{\Delta u}_{L^2}^2 + \alpha^2\norm{\nabla\Delta u}_{L^2}^2 + \norm{\Delta b}_{L^2}^2\big)
   + 
   \nu \norm{\nabla\Delta u}_{L^2}^2 
   +
   \eta \norm{\nabla\Delta b}_{L^2}^2 
   \\&
   \leq
   C\big(\norm{\Delta u}_{L^2}^2 + \alpha^2\norm{\nabla\Delta u}_{L^2}^2 + \norm{\Delta b}_{L^2}^2\big),
\end{align*}
where the above constant $C$ depends on $\nu, \eta$, and the $H^2$-bounds on $u$ and $b$. Hence, Gr\"onwall's inequality leads to the desired conclusion, i.e., we have proved $$u\in C(0, T; H^3\cap V).\quad\text{\qedsymbol}$$

Next, we prove the $H^3$-bound on $b$, and to do so, we multiply $\partial^2 u$ to \eqref{VVV_MHD_u} after applying the operator $\partial^2$ to it first, and multiply $\partial^3 b$ to \eqref{VVV_MHD_b} after applying the operator $\partial^3$ to it first, and then, we add the two resulting equations, integrate by parts, so that we get
\begin{align*}
    &
    \frac{1}{2} \frac{d}{dt}\big(\norm{\partial^{2} u}_{L^2}^2 + \alpha^2\norm{\nabla\partial^{2} u}_{L^2}^2 +  \norm{\partial^3 b}_{L^2}^2\big)
    + 
    \nu \norm{\nabla\partial^{2} u}_{L^2}^2 
    +
    \eta \norm{\nabla\partial^3 b}_{L^2}^2 
    \\&
    = 
    -\int_{\mathbb{T}^3}\big(\partial^{2}(w\times u)\big)\cdot(\partial{^2} u) 
    + 
    \int_{\mathbb{T}^3}\big(\partial^{2}(b\cdot \nabla)b\big)\cdot(\partial{^2} u)
    \\&\quad
    +
    \int_{\mathbb{T}^3}\partial^3\big((b \cdot \nabla) u\big)\cdot (\partial^3 b) 
    - 
    \int_{\mathbb{T}^3}\partial^3\big((u \cdot \nabla) b\big)\cdot (\partial^3 b).
\end{align*}
We then estimate each of the four integrals on the right side of the above equations. Note that after integration by parts, the first two integrals are bounded the same way as in the $H^3$-estimates of $u$, so for simplicity, we focus on the last two integrals. For the third integral, we first apply the Leibniz rule to rewrite it as  
\begin{align*}
    &
    \int_{\mathbb{T}^3}\partial^3\big((b \cdot \nabla) u\big)\cdot (\partial^3 b)
    =
    \sum\limits_{0\leq l\leq 3}\binom{3}{l}\int_{\mathbb{T}^3}\big((\partial^{l}b \cdot \nabla) \partial^{3-l}u\big)\cdot (\partial^3 b)
    \\&
    =
    {\rm I}_0
    + 
    {\rm I}_1
    + 
    {\rm I}_2
    +
    {\rm I}_3.
\end{align*}
Then, we estimate each ${\rm I}_{j}$, $j=1,2,3,4$ as follows. For ${\rm I}_0$, we integrate by parts and bounded it by
\begin{align*}
    |{\rm I}_0|
    &
    \leq
    C\int_{\mathbb{T}^3}|b| |\partial^{3} u| |\nabla\partial^{3} b|
    \leq
    C\normLp{\infty}{b}\normL{\partial^{3} u}\normL{\nabla\partial^{3} b}
    \\&
    \leq
    C\normL{\nabla b}^{1/2}\normL{\Delta b}^{1/2}\big(\normL{\partial^{3} u}^{1/2}\big)^2\normL{\nabla\partial^{3} b}
    \\&
    \leq
    C\normL{\nabla b}\normL{\Delta b}\normL{\partial^{3} u}
    +
    \frac{\eta}{8}\normL{\nabla\partial^{3} b}^2
    \\&
    \leq
    C\normL{\nabla b}^2\normL{\Delta b}^2
    +
    \frac{\nu}{8}\normL{\nabla\partial^{2} u}^2
    +
    \frac{\eta}{8}\normL{\nabla\partial^{3} b}^2
    \\&
    \leq
    C\normL{\partial^{3} b}^2
    +
    \frac{\nu}{8}\normL{\nabla\partial^{2} u}^2
    +
    \frac{\eta}{8}\normL{\nabla\partial^{3} b}^2,
\end{align*}
where we used the boundedness of ${\norm{b}}_{H^1}$ and ${\norm{u}}_{H^3}$, as well as Sobolev, Poincar\'e's, and Young's inequalities. The integral ${\rm I}_1$ is bounded by
\begin{align*}
    |{\rm I}_1|
    &
    \leq
    C\int_{\mathbb{T}^3}|\nabla b| |\nabla\partial^{2} u| |\partial^{3} b|
    \leq
    C\normLp{3}{\nabla b}\normL{\partial^{3} u}\normLp{6}{\partial^{3} b}
    \\&
    \leq
    C\normL{\partial^{3} b}^2
    +
    \frac{\nu}{8}\normL{\nabla\partial^{2} u}^2
    +
    \frac{\eta}{8}\normL{\nabla\partial^{3} b}^2,
\end{align*}
where the detailed rationale of the last inequality is similar to that of the estimates on ${\rm I}_0$. Regarding ${\rm I}_2$, we estimate as
\begin{align*}
    |{\rm I}_2|
    &
    \leq
    C\int_{\mathbb{T}^3}|\nabla^{2} b| |\partial^{2} u| |\partial^{3} b|
    \leq
    C\normL{\partial^{2} b}\normLp{3}{\partial^{2} u}\normLp{6}{\partial^{3} b}
    \\&
    \leq
    C{\norm{b}}_{H^2}\normL{\partial^{2} u}^{1/2}\normL{\nabla\partial^{2} u}^{1/2}\normL{\nabla\partial^{3} b}
    \\&
    \leq
    C{\norm{b}}_{H^2}^2\normL{\partial^{2} u}\normL{\nabla\partial^{2} u}
    +
    \frac{\eta}{8}\normL{\nabla\partial^{3} b}^2
    \\&
    \leq
    C\normL{\partial^{2} u}^2
    +
    \frac{\nu}{8}\normL{\nabla\partial^{2} u}^2
    +
    \frac{\eta}{8}\normL{\nabla\partial^{3} b}^2.
\end{align*}
As for the integral ${\rm I}_3$, we proceed analogously as 
\begin{align*}
    |{\rm I}_3|
    &
    \leq
    C\int_{\mathbb{T}^3}|\nabla^{3} b| |\nabla u| |\partial^{3} b|
    \leq
    C\normL{\partial^{3} b}\normLp{\infty}{\nabla u}\normL{\partial^{3} b}
    \\&
    \leq
    C\normL{\nabla u}^{1/2}\normL{\Delta u}^{1/2}\normL{\partial^{3} b}\normL{\partial^{3} b}
    \\&
    \leq
    C\normL{\partial^{3} b}^2
    +
    \frac{\eta}{8}\normL{\nabla\partial^{3} b}^2,
\end{align*}
where we used the boundedness of ${\norm{u}}_{H^3}$. 

As for the last integral, we proceed similarly and rewrite it as 
\begin{align*}
    &
    \int_{\mathbb{T}^3}\partial^3\big((b \cdot \nabla) u\big)\cdot (\partial^3 b)
    =
    \sum\limits_{1\leq l\leq 3}\binom{3}{l}\int_{\mathbb{T}^3}\big((\partial^{l}u \cdot \nabla) \partial^{3-l}b\big)\cdot (\partial^3 b)
    \\&
    =
    {\rm II}_1
    +
    {\rm II}_2
    +
    {\rm II}_3,
\end{align*}
where we note that the term with $l=0$ vanishes due to the divergence-free condition on $b$. Then, we estimate each ${\rm II}_{j}$, $j=1,2,3$ as follows. 
For ${\rm II}_1$, we have 
\begin{align*}
    |{\rm II}_1|
    &
    \leq
    C\int_{\mathbb{T}^3}|\nabla u| |\nabla\partial^{2} b| |\partial^{3} b|
    \leq
    C\normLp{\infty}{\nabla u}\normL{\partial^{3} b}^2
    \\&
    \leq
    C{\norm{u}}_{H^3}\normL{\partial^{3} b}^2
    \leq
    C\normL{\partial^{3} b}^2. 
\end{align*}
Next, ${\rm II}_2$ is bounded by 
\begin{align*}
    |{\rm II}_2|
    &
    \leq
    C\int_{\mathbb{T}^3}|\partial^{2} u| |\partial^{2} b| |\partial^{3} b|
    \leq
    C\normL{\partial^{2} u}\normLp{3}{\partial^{2} b}\normLp{6}{\partial^{3} b}
    \\&
    \leq
    C\normL{\partial^{2} u}\normL{\partial^{3} b}\normL{\nabla\partial^{3} b}
    \\&
    \leq
    C\normL{\partial^{3} b}^2
    +
    \frac{\eta}{8}\normL{\nabla\partial^{3} b}^2, 
\end{align*}
where we used the boundedness of ${\norm{u}}_{H^2}$. Lastly, we estimate ${\rm II}_3$ as 
\begin{align*}
    |{\rm II}_3|
    &
    \leq
    C\int_{\mathbb{T}^3}|\partial^{3} u| |\nabla b| |\partial^{3} b|
    \leq
    C\normL{\partial^{3} u}\normLp{6}{\nabla b}\normLp{3}{\partial^{3} b}
    \\&
    \leq
    C\normL{\Delta b}\normL{\partial^{3} b}^{1/2}\normL{\nabla\partial^{3} b}^{1/2}\normL{\nabla\partial^{2} u}
    \\&
    \leq
    C{\norm{b}}_{H^2}^2\normL{\partial^{3} b}\normL{\nabla\partial^{3} b}
    +
    \frac{\nu}{8}\normL{\nabla\partial^{2} u}^2
    \\&
    \leq
    C\normL{\partial^{3} b}^2
    +
    \frac{\eta}{8}\normL{\nabla\partial^{3} b}
    +
    \frac{\nu}{8}\normL{\nabla\partial^{2} u}^2. 
\end{align*}
Therefore, by combining all the above bounds, we obtain
\begin{align*}
   &\frac{d}{dt}\big(\norm{\partial^{2} u}_{L^2}^2 + \alpha^2\norm{\nabla\partial^{2} u}_{L^2}^2 + \norm{\partial^{3} b}_{L^2}^2\big)
   + 
   \nu\norm{\nabla\partial^{2} u}_{L^2}^2 
   +
   \eta\norm{\nabla\partial^{3} b}_{L^2}^2 
   \\&
   \leq
   C\big(\norm{\partial^{2} u}_{L^2}^2 + \alpha^2\norm{\nabla\partial^{2} u}_{L^2}^2 + \norm{\partial^{3} b}_{L^2}^2\big),
\end{align*}
where the above constant $C$ depends on $\nu, \eta$, and the $H^3$-bounds on $u$ and $b$. Hence, Gr\"onwall's inequality implies the desired conclusion. Namely, we have proved $$b\in C(0, T; H^3\cap V)\cap L^2(0, T; H^4\cap V).\quad\text{\qedsymbol}$$

\subsection{\texorpdfstring{$H^2$}--estimates of \texorpdfstring{$w$}-}\label{w-H2}
Now, we are ready to obtain the bound on the $H^2$-norm of $w$, and to do so, we multiple equations \eqref{VVV_MHD_u}, \eqref{VVV_MHD_w}, and \eqref{VVV_MHD_b} by $\Delta^2 u$ and $\Delta^2 w$, and $\Delta^2 b$, respectively, integrate by parts, and add the three resulted equations together, so that we obtain 
\begin{align*}
    &
    \frac{1}{2} \frac{d}{dt}\big(\norm{\Delta u}_{L^2}^2 + \alpha^2\norm{\nabla\Delta u}_{L^2}^2 + \norm{\Delta w}_{L^2}^2 + \norm{\Delta b}_{L^2}^2\big)
    \\&\quad
    + 
    \nu \norm{\nabla\Delta u}_{L^2}^2 
    +
    \nu \norm{\nabla\Delta w}_{L^2}^2 
    +
    \eta \norm{\nabla\Delta b}_{L^2}^2 
    \\&
    = 
    -\int_{\mathbb{T}^3}(w\times u)\cdot(\Delta^2 u) 
    + 
    \int_{\mathbb{T}^3}(b\cdot \nabla)b\cdot(\Delta^2 u)
    \\&\quad
    -
    \int_{\mathbb{T}^3}(u\cdot\nabla)w\cdot(\Delta^2 w)
    +
    \int_{\mathbb{T}^3}(w\cdot\nabla)u\cdot(\Delta^2 w)
    +
    \int_{\mathbb{T}^3}\nabla\times(j\times b)\cdot(\Delta^2 w)
    \\&\quad\quad 
    - 
    \int_{\mathbb{T}^3}(u \cdot \nabla) b\cdot (\Delta^2 b)
    +
    \int_{\mathbb{T}^3}(b \cdot \nabla) u\cdot (\Delta^2 b).
\end{align*}
Note that the second and the last two integrals on the right side of the above equation are estimated similarly as in the $H^3$-estimates of $u$ and the $H^2$-estimates of $b$, so we integrate by parts and bound the first integral by 
\begin{align*}
    &
    \int_{\mathbb{T}^3}|\nabla w||u||\nabla\Delta u|
    +
    \int_{\mathbb{T}^3} |w||\nabla u||\nabla\Delta u|
    \\&
    \leq
    C\normLp{3}{\nabla w}\normLp{6}{u}\normL{\nabla\Delta u}
    +
    C\normLp{\infty}{w}\normL{\nabla u}\normL{\nabla\Delta u}
    \\&
    \leq
    C\normL{\nabla w}^{1/2}\normL{\Delta w}^{1/2}\normL{\nabla u}\normL{\nabla\Delta u}
    +
    C\normL{\Delta w}\normL{\nabla u}\normL{\nabla\Delta u}
    \\&
    \leq
    C\normL{\Delta w}^2
    +
    \frac{\nu}{8}\normL{\nabla\Delta u}^2,
\end{align*} 
where we used Sobolev, Poincar\,e's, and Young's inequalities, as well as the boundedness of the $H^1$-norm of $u$. 

Next, we also integrate by parts before bounding the third integral by
\begin{align*}
    &
    \int_{\mathbb{T}^3}|\nabla w||\nabla u||\nabla\Delta w|
    +
    \int_{\mathbb{T}^3} |w||\nabla^2 u||\nabla\Delta u|
    \\&
    \leq
    C\normLp{3}{\nabla w}\normLp{6}{\nabla u}\normL{\nabla\Delta u}
    +
    C\normLp{\infty}{w}\normL{\nabla u}\normL{\nabla\Delta u}
    \\&
    \leq
    C\norm{\Delta u}_{H^2}\normL{\Delta w}\normL{\nabla\Delta u}
    \\&
    \leq
    C\normL{\Delta w}^2
    +
    \frac{\nu}{8}\normL{\nabla\Delta u}^2,
\end{align*}
where we also used the boundedness of $\norm{u}_{H^2}$. 

As for the fourth integral, we proceed similarly and estimate it as 
\begin{align*}
    &
    \int_{\mathbb{T}^3}|\nabla u||\nabla w||\nabla\Delta w|
    +
    \int_{\mathbb{T}^3} |u||\nabla^2 w||\nabla\Delta u|
    \\&
    \leq
    C\normLp{6}{\nabla u}\normLp{3}{\nabla w}\normL{\nabla\Delta u}
    +
    C\normLp{\infty}{u}\normL{\nabla w}\normL{\nabla\Delta u}
    \\&
    \leq
    C\norm{\Delta u}_{H^2}\normL{\Delta w}\normL{\nabla\Delta u}
    \\&
    \leq
    C\normL{\Delta w}^2
    +
    \frac{\nu}{8}\normL{\nabla\Delta u}^2.
\end{align*}
As for the fifth integral, we integrate by parts, apply Leibniz rule, as well as the fact that the fact $\nabla\cdot b=0$ together with the vector calculus identity 
\begin{equation}\label{calculus_identity}
\nabla\times(\vec{F}\times\vec{G}) = \vec{F}(\nabla\cdot\vec{G}) - \vec{G}(\nabla\cdot\vec{F}) + (\vec{G}\cdot\nabla)\vec{F} - (\vec{F}\cdot\nabla)\vec{G},
\end{equation}
where $\vec{F}$ and $\vec{G}$ are any twice-differentiable vector fields, so that we have
\begin{align*}
    &
    \Big|\int_{\mathbb{T}^3}\nabla\times(j\times b)\cdot(\Delta^2 w)\Big|
    =
    \Big|\int_{\mathbb{T}^3}\nabla\big(\nabla\times(j\times b)\big)\cdot(\nabla\Delta w)\Big|
    \\&
    \leq
    \int_{\mathbb{T}^3}|\nabla b||\nabla j||\nabla\Delta w|
    +
    \int_{\mathbb{T}^3}|b||\nabla^2 j||\nabla\Delta w|
    +
    \int_{\mathbb{T}^3}|\nabla^2 b||j||\nabla\Delta w|
    \\&
    \leq
    C\normLp{6}{\nabla b}\normLp{3}{\nabla j}\normL{\nabla\Delta w}
    +
    C\normLp{\infty}{b}\normL{\nabla^2 j}\normL{\nabla\Delta w}
    \\&\quad
    +
    C\normLp{3}{\nabla^2 b}\normLp{6}{j}\normL{\nabla\Delta w}
    \\&
    \leq
    C\norm{b}_{H^3}\normL{\Delta b}\normL{\nabla\Delta w}
    \leq
    C\normL{\Delta b}^2
    +
    \frac{\nu}{8}\normL{\nabla\Delta w}^2,
\end{align*}
where we also used the $H^3$-boundedness of $b$. 

Thus, combining all the above estimates, we have  
\begin{align*}
   &\frac{d}{dt}\big(\norm{\Delta u}_{L^2}^2 + \alpha^2\norm{\nabla\Delta u}_{L^2}^2 + \norm{\Delta w}^2 + \norm{\Delta b}_{L^2}^2\big)
   \\&\quad
   + 
   \nu \norm{\nabla\Delta u}_{L^2}^2 
   +
   \nu \norm{\nabla\Delta w}_{L^2}^2
   +
   \eta \norm{\nabla\Delta b}_{L^2}^2 
   \\&
   \leq
   C\big(\norm{\Delta u}_{L^2}^2 + \alpha^2\norm{\nabla\Delta u}_{L^2}^2 + \norm{\Delta w}^2 + \norm{\Delta b}_{L^2}^2\big),
\end{align*}
where the above constant $C$ depends on $\nu, \eta$, and the $H^3$-bounds on $u$ and $b$. Hence, Gr\"onwall's inequality implies the desired conclusion, i.e., we have proved $$w\in C(0, T; H^2\cap V)\cap L^2(0, T; H^3\cap V).$$

Following the above bootstrapping strategy, i.e., by first proving the higher-order-norm estimates on $u$, then on $b$, we are able to obtain the higher-order estimates on $w$ in $H^s$ for $s=3, 4, \dots$. Thus, Theorem~\ref{thm-regularity} is now proved. 
\qedsymbol

\subsection{Proof of Uniqueness}\label{uniqueness}
Following the well-known weak-strong uniqueness argument for the Navier-Stokes equations, it suffices to show that the strong solution is unique for all $\alpha>0$. Suppose there exist two solutions $(u,w,b,p)$ and $(u',w',b',p')$ to System~\eqref{VVV_MHD}, with the same initial data $u_0, b_0\in D(A)$ and $w_0 \in V$, and with $j=\nabla\times b$ and $j'=\nabla\times b'$, on their common time-interval of existence. By denoting the difference $\tilde{u}=u-u', \tilde{w}=w-w', \tilde{b}=b-b'$, and $\tilde{p}=p-p'$, and subtracting the two sets of equations satisfied by the two solutions, we obtain 
\begin{equation}
     \left\{
     \begin{aligned}
       &
       (I - \alpha^2\Delta)\frac{\partial \tilde{u}}{\partial t} 
       - 
       \nu\Delta \tilde{u} 
       + 
       \tilde{w}\times u 
       +
       w'\times \tilde{u}
       + 
       \nabla \tilde{p} 
       = 
       \tilde{j}\times b
       + 
       j'\times \tilde{b},
       \\&
       \frac{\partial \tilde{w}}{\partial t} 
       - 
       \nu\Delta \tilde{w} 
       + 
       (u\cdot\nabla) \tilde{w}
       + 
       (\tilde{u}\cdot\nabla) w'
       - 
       (w\cdot\nabla) \tilde{u}
       - 
       (\tilde{w}\cdot\nabla)u'
       \\&\quad
       = 
       \nabla\times(\tilde{j}\times b)
       +
       \nabla\times (j'\times \tilde{b}),
       \\&
       \frac{\partial \tilde{b}}{\partial t} 
       - 
       \eta\Delta \tilde{b} 
       + 
       (u\cdot\nabla) \tilde{b}
       + 
       (\tilde{u}\cdot\nabla)b'
       = 
       (b\cdot\nabla) \tilde{u}
       +
       (\tilde{b}\cdot\nabla)u', 
       \\&
       \nabla \cdot u' = \nabla \cdot w' = \nabla \cdot b' = 0,
    \end{aligned}
    \right.
    \label{Sys_Unique}
  \end{equation}
where $\tilde{j}=\nabla\times\tilde{b}$, and we used the alternative form of the $u$-equation in System~\eqref{VVV_MHD} corresponding to \eqref{curl_form} in System~\eqref{MHD}, i.e., with $j\times b$ in place of $(b\cdot \nabla) b$. 

Next, we multiply each of the equations of $\tilde{u}, \tilde{w}$ and $\tilde{b}$ by $\tilde{u}, \tilde{w}$, $\tilde{b}$, respectively, and also multiply the $\tilde{u}$-equation by $A\tilde{u}=-\Delta\tilde{u}$ and the $\tilde{b}$-equation by $A\tilde{b}=-\Delta\tilde{b}$; then, we integrate by parts on $\mathbb{T}^3$, and add all five equations, so we obtain 
\begin{align*}
  &
  \frac{1}{2}\frac{d}{dt}\big(\|\tilde{u}\|_{L^2}^2 + \alpha^2\|\nabla\tilde{u}\|_{L^2}^2 + \|\tilde{w}\|_{L^2}^2 + \|\tilde{b}\|_{L^2}^2 + \|\nabla \tilde{b}\|_{L^2}^2\big)
  \\&\quad
  +
  \nu\|\nabla\tilde{u}\|_{L^2}^2
  +
  \nu\|\nabla\tilde{w}\|_{L^2}^2
  +
  \eta\|\nabla\tilde{b}\|_{L^2}^2
  +
  \eta\|\Delta\tilde{b}\|_{L^2}^2
  \\&
  =
  -
  \int_{\mathbb{T}^3}(\tilde{w}\times u)\cdot\tilde{u}
  +
  \int_{\mathbb{T}^3}(\tilde{j}\times b)\cdot\tilde{u}
  +
  \int_{\mathbb{T}^3}(j'\times \tilde{b})\cdot\tilde{u}
  \\&\quad
  -
  \int_{\mathbb{T}^3}(\tilde{u}\cdot \nabla) w'\cdot\tilde{w}
  +
  \int_{\mathbb{T}^3}(w\cdot \nabla) \tilde{u}\cdot\tilde{w}
  +
  \int_{\mathbb{T}^3}(\tilde{w}\cdot \nabla) u'\cdot\tilde{w}
  \\&\quad
  +
  \int_{\mathbb{T}^3}(b\cdot \nabla) \tilde{j}\cdot\tilde{w}
  -
  \int_{\mathbb{T}^3}(\tilde{j}\cdot \nabla) b\cdot\tilde{w}
  +
  \int_{\mathbb{T}^3}(\tilde{b}\cdot \nabla) j'\cdot\tilde{w}
  -
  \int_{\mathbb{T}^3}(j'\cdot \nabla) \tilde{b}\cdot\tilde{w}
  \\&\quad
  -
  \int_{\mathbb{T}^3}(\tilde{u}\cdot \nabla) b'\cdot\tilde{b}
  +
  \int_{\mathbb{T}^3}(b\cdot \nabla) \tilde{u}\cdot\tilde{b}
  +
  \int_{\mathbb{T}^3}(\tilde{b}\cdot \nabla) u'\cdot\tilde{b}
  \\&\quad
  +
  \int_{\mathbb{T}^3}(\tilde{w}\times u)\cdot\Delta\tilde{u}
  +
  \int_{\mathbb{T}^3}(w'\times\tilde{u})\cdot\Delta\tilde{u}
  -
  \int_{\mathbb{T}^3}(\tilde{j}\times b)\cdot\Delta\tilde{u}
  -
  \int_{\mathbb{T}^3}(j'\times \tilde{b})\cdot\Delta\tilde{u}
  \\&\quad
  +
  \int_{\mathbb{T}^3}(u\cdot \nabla) \tilde{b}\cdot\Delta\tilde{b}
  +
  \int_{\mathbb{T}^3}(\tilde{u}\cdot \nabla) b'\cdot\Delta\tilde{b}
  -
  \int_{\mathbb{T}^3}(b\cdot \nabla) \tilde{u}\cdot\Delta\tilde{b}
  -
  \int_{\mathbb{T}^3}(\tilde{b}\cdot \nabla) u'\cdot\Delta\tilde{b}
  \\&
  :=\Lambda_1 + \Lambda_2 + \cdots + \Lambda_{21},
\end{align*}
where we used the vector calculus identity (\eqref{calculus_identity})
and the fact $\nabla\cdot(\nabla\times \vec{F})=0$ for any twice-differentiable vector fields $\vec{F}$. 

Then, we estimate the seventeen terms $\Lambda_1$ through $\Lambda_{17}$ on the right side of the above equation. 
First, for $\Lambda_1$, by H\"older and Cauchy-Schwarz inequalities, we have 
\begin{align*}
   |\Lambda_1| 
   &
   \leq
   \|u\|_{L^2}\|\tilde{u}\|_{L^3}\|\tilde{w}\|_{L^6}
   \leq
   C\|u\|_{L^2}\|\tilde{u}\|_{L^2}^{1/2}\|\nabla\tilde{u}\|_{L^2}^{1/2}\|\nabla\tilde{w}\|_{L^2}
   \\&
   \leq
   C\|u\|_{L^2}^4\|\tilde{u}\|_{L^2}^2
   +
   \frac{\nu}{16}\|\nabla\tilde{u}\|_{L^2}^2
   +
   \frac{\nu}{16}\|\nabla\tilde{w}\|_{L^2}^2,
\end{align*}
where we alse used Sobolev inequality. 
Similarly, $\Lambda_2$ is estimated as 
\begin{align*}
   |\Lambda_2| 
   &
   \leq
   \|\tilde{j}\|_{L^2}\|b\|_{L^6}\|\tilde{u}\|_{L^3}
   \leq
   C\|\nabla\tilde{b}\|_{L^2}\|\nabla b\|_{L^2}\|\|\tilde{u}\|_{L^2}^{1/2}\|\nabla\tilde{u}\|_{L^2}^{1/2}
   \\&
   \leq
   C\|\nabla b\|_{L^2}^2 \|\tilde{u}\|_{L^2}\|\nabla\tilde{u}\|_{L^2}
   +
   \frac{\eta}{16}\|\nabla\tilde{b}\|_{L^2}^2
   \\&
   \leq
   C\|\nabla b\|_{L^2}^4\|\tilde{u}\|_{L^2}^2
   +
   \frac{\nu}{16}\|\nabla\tilde{u}\|_{L^2}^2
   +
   \frac{\eta}{16}\|\nabla\tilde{b}\|_{L^2}^2.
\end{align*}
Estimates on $\Lambda_3$ is also analogous to $\Lambda_1$ and $\Lambda_2$, and we have
\begin{align*}
   |\Lambda_3| 
   &
   \leq
   C\|\nabla b'\|_{L^2}^4\|\tilde{b}\|_{L^2}^2
   +
   \frac{\eta}{16}\|\nabla\tilde{b}\|_{L^2}^2
   +
   \frac{\nu}{16}\|\nabla\tilde{u}\|_{L^2}^2.
\end{align*}
Regarding $\Lambda_4$, we also use similar arguments as above to obtain 
\begin{align*}
   |\Lambda_4|
   &
   \leq
   C\|\nabla w'\|_{L^2}\|\tilde{u}\|_{L^3}\|\tilde{w}\|_{L^6}
   \leq
   C\|\nabla w'\|_{L^2}\|\tilde{u}\|_{L^2}^{1/2}\|\nabla\tilde{u}\|_{L^2}^{1/2}\|\nabla\tilde{w}\|_{L^2}
   \\&
   \leq
   C\|\nabla w'\|_{L^2}^2\|\tilde{u}\|_{L^2}\|\nabla\tilde{u}\|_{L^2}
   +
   \frac{\nu}{16}\|\nabla\tilde{w}\|_{L^2}^2
   \\&
   \leq
   C\|\nabla w'\|_{L^2}^4 \|\tilde{u}\|_{L^2}^2 
   +
   \frac{\nu}{16}\|\nabla\tilde{u}\|_{L^2}^2
   +
   \frac{\nu}{16}\|\nabla\tilde{w}\|_{L^2}^2, 
\end{align*}
while estimates on $\Lambda_5$ proceed similarly as those on $\Lambda_4$, i.e., we have 
\begin{align*}
   |\Lambda_5|
   &
   \leq
   C\|\nabla w\|_{L^2}^4\|\tilde{w}\|_{L^2}^2
   +
   \frac{\nu}{16}\|\nabla\tilde{u}\|_{L^2}^2
   +
   \frac{\nu}{16}\|\nabla\tilde{w}\|_{L^2}^2. 
\end{align*}
Then, we estimate $\Lambda_6$ as 
\begin{align*}
   |\Lambda_6|
   &
   \leq
   \|\tilde{w}\|_{L^3}\|\nabla u'\|_{L^2}\|\tilde{w}\|_{L^6}
   \leq
   C\|\nabla u'\|_{L^2}\|\tilde{w}\|_{L^2}^{1/2}\|\nabla\tilde{w}\|_{L^2}^{3/2}
   \\&
   \leq
   C\|\nabla u'\|_{L^2}^4\|\tilde{w}\|_{L^2}^2
   +
   \frac{\nu}{16}\|\nabla\tilde{w}\|_{L^2}^2,
\end{align*}
where we have frequently used Young's inequality. 

Next, we bound each of the next four terms as  
\begin{align*}
   |\Lambda_{7}| 
   &
   \leq
   \|b\|_{L^6}\|\nabla\tilde{j}\|_{L^2}\|\tilde{w}\|_{L^3}
   \leq
   C\|\nabla b\|_{L^2} \|\Delta\tilde{b}\|_{L^2} \|\tilde{w}\|_{L^2}^{1/2}\|\nabla\tilde{w}\|_{L^2}^{1/2}
   \\&
   \leq
   C\|\nabla b\|_{L^2}^2 \|\tilde{w}\|_{L^2}\|\nabla\tilde{w}\|_{L^2}
   +
   \frac{\eta}{8}\|\Delta\tilde{b}\|_{L^2}^2
   \\&
   \leq
   C\|\nabla b\|_{L^2}^4 \|\tilde{w}\|_{L^2}^2
   +
   \frac{\nu}{16}\|\nabla\tilde{w}\|_{L^2}^2
   +
   \frac{\eta}{8}\|\Delta\tilde{b}\|_{L^2}^2, 
\end{align*}
and similarly, 
\begin{align*}
   |\Lambda_{8}| 
   &
   \leq
   \|\tilde{j}\|_{L^6}\|\nabla b\|_{L^2}\|\tilde{w}\|_{L^3}
   \\&
   \leq
   C\|\nabla b\|_{L^2}^4 \|\tilde{w}\|_{L^2}^2
   +
   \frac{\eta}{8}\|\Delta\tilde{b}\|_{L^2}^2
   +
   \frac{\nu}{16}\|\nabla\tilde{w}\|_{L^2}^2,
\end{align*}
while  
\begin{align*}
   |\Lambda_{9} | 
   &
   \leq
   \|\tilde{b}\|_{L^3}\|\nabla j'\|_{L^2}\|\tilde{w}\|_{L^6}
   \leq
   C\|\Delta b'\|_{L^2}\|\tilde{b}\|_{L^2}^{1/2}\|\nabla\tilde{b}\|_{L^2}^{1/2}\|\nabla\tilde{w}\|_{L^2}
   \\&
   \leq
   C\|\Delta b'\|_{L^2}^4\|\tilde{b}\|_{L^2}^2
   +
   \frac{\nu}{16}\|\nabla\tilde{w}\|_{L^2}^2
   +
   \frac{\eta}{16}\|\nabla\tilde{b}\|_{L^2}^2,
\end{align*}
and
\begin{align*}
   |\Lambda_{10} | 
   &
   \leq
   \|\nabla b'\|_{L^6}\|\nabla\tilde{b}\|_{L^2}\|\tilde{w}\|_{L^3}
   \leq
   C\|\Delta b'\|_{L^2}\|\tilde{w}\|_{L^2}^{1/2}\|\nabla\tilde{w}\|_{L^2}^{1/2}\|\nabla\tilde{b}\|_{L^2}
   \\&
   \leq
   C\|\Delta b'\|_{L^2}^4\|\tilde{w}\|_{L^2}^2
   +
   \frac{\eta}{16}\|\nabla\tilde{b}\|_{L^2}^2
   +
   \frac{\nu}{16}\|\nabla\tilde{w}\|_{L^2}^2.
\end{align*}

Estimates on the next three terms $\Lambda_{11}, \Lambda_{12}$ and $\Lambda_{13}$ proceed analogously to those of the first six terms, so we have 
\begin{align*}
   |\Lambda_{11}| 
   &
   \leq
   C\|\nabla b'\|_{L^2}^4\|\tilde{u}\|_{L^2}^2
   +
   \frac{\eta}{16}\|\nabla\tilde{b}\|_{L^2}^2
   +
   \frac{\nu}{16}\|\nabla\tilde{u}\|_{L^2}^2,
\end{align*}
and
\begin{align*}
   |\Lambda_{12}| 
   &
   \leq
   C\|\nabla b\|_{L^2}^4\|\tilde{b}\|_{L^2}^2
   +
   \frac{\nu}{16}\|\nabla\tilde{u}\|_{L^2}^2
   +
   \frac{\eta}{16}\|\nabla\tilde{b}\|_{L^2}^2,
\end{align*}
and 
\begin{align*}
   |\Lambda_{13} | 
   &
   \leq
   C\|\nabla u'\|_{L^2}^4\|\tilde{b}\|_{L^2}^2
   +
   \frac{\eta}{16}\|\nabla\tilde{b}\|_{L^2}^2.
\end{align*}
We then estimate terms $\Lambda_{14}$ through $\Lambda_{17}$ as follows. 
\begin{align*}
   |\Lambda_{14}| 
   &
   \leq
   \|\tilde{w}\|_{L^3}\|u\|_{L^6}\|\Delta\tilde{u}\|_{L^2}
   \leq
   C\|\nabla u\|_{L^2}\|\tilde{w}\|_{L^2}^{1/2}\|\nabla\tilde{w}\|_{L^2}^{1/2}\|\Delta\tilde{u}\|_{L^2}
   \\&
   \leq
   C\|\nabla u\|_{L^2}^4\|\tilde{w}\|_{L^2}^2
   +
   \frac{\nu}{8}\|\Delta\tilde{u}\|_{L^2}^2
   +
   \frac{\nu}{16}\|\nabla\tilde{w}\|_{L^2}^2,
\end{align*}
and $\Lambda_{15}$ is bounded similarly by
\begin{align*}
   |\Lambda_{15}| 
   &
   \leq
   C\|\nabla w'\|_{L^2}^4\|\tilde{u}\|_{L^2}^2
   +
   \frac{\nu}{8}\|\Delta\tilde{u}\|_{L^2}^2
   +
   \frac{\nu}{16}\|\nabla\tilde{u}\|_{L^2}^2,
\end{align*}
while 
\begin{align*}
   |\Lambda_{16}| 
   &
   \leq
   \|\tilde{j}\|_{L^3}\|b\|_{L^6}\|\Delta\tilde{u}\|_{L^2}
   \leq
   C\|\nabla b\|_{L^2}\|\nabla\tilde{b}\|_{L^2}^{1/2}\|\Delta\tilde{b}\|_{L^2}^{1/2}\|\Delta\tilde{u}\|_{L^2}
   \\&
   \leq
   C\|\nabla b\|_{L^2}^4\|\nabla\tilde{b}\|_{L^2}^2
   +
   \frac{\nu}{8}\|\Delta\tilde{u}\|_{L^2}^2
   +
   \frac{\eta}{8}\|\Delta\tilde{b}\|_{L^2}^2, 
\end{align*}
while 
\begin{align*}
   |\Lambda_{17}| 
   &
   \leq
   \|j'\|_{L^6}\|\tilde{b}\|_{L^3}\|\Delta\tilde{u}\|_{L^2}
   \leq
   C\|\Delta b'\|_{L^2}\|\tilde{b}\|_{L^2}^{1/2}\|\nabla\tilde{b}\|_{L^2}^{1/2}\|\Delta\tilde{u}\|_{L^2}
   \\&
   \leq
   C\|\Delta b'\|_{L^2}^4\|\tilde{b}\|_{L^2}^2
   +
   \frac{\nu}{8}\|\Delta\tilde{u}\|_{L^2}^2
   +
   \frac{\eta}{16}\|\nabla\tilde{b}\|_{L^2}^2. 
\end{align*}
Finally, we follow similar reasoning as well as Agmon's inequality and the boundedness of the $H^2$-norms of $u'$ and $b'$, and estimate the last four terms $\Lambda_{18}$ through $\Lambda_{21}$ as 
\begin{align*}
   |\Lambda_{18}| 
   &
   \leq
   \|u\|_{L^6}\|\nabla\tilde{b}\|_{L^3}\|\Delta\tilde{b}\|_{L^2}
   \leq
   C\|\nabla u\|_{L^2}\|\nabla\tilde{b}\|_{L^2}^{1/2}\|\Delta\tilde{b}\|_{L^2}^{3/2}
   \\&
   \leq
   C\|\nabla u\|_{L^2}^4\|\nabla\tilde{b}\|_{L^2}^2
   +
   \frac{\eta}{16}\|\Delta\tilde{b}\|_{L^2}^2,
\end{align*}
and
\begin{align*}
   |\Lambda_{19}| 
   &
   \leq
   \|\tilde{u}\|_{L^3}\|\nabla b'\|_{L^6}\|\Delta\tilde{b}\|_{L^2}
   \leq
   C \|\Delta b'\|_{L^2}\|\tilde{u}\|_{L^2}^{1/2}\|\nabla\tilde{u}\|_{L^2}^{1/2}\|\Delta\tilde{b}\|_{L^2}
   \\&
   \leq
   C\|\Delta b'\|_{L^2}^4 \|\tilde{u}\|_{L^2}^2
   +
   \frac{\nu}{16}\|\nabla\tilde{u}\|_{L^2}^2
   +
   \frac{\eta}{8}\|\Delta\tilde{b}\|_{L^2}^2,
\end{align*}
and 
\begin{align*}
   |\Lambda_{20} | 
   &
   \leq
   \|b\|_{L^{6}}\|\nabla\tilde{u}\|_{L^3}\|\Delta\tilde{b}\|_{L^2}
   \leq
   C\|\nabla b\|_{L^2}^4\|\nabla\tilde{u}\|_{L^2}^2
   +
   \frac{\nu}{8}\|\Delta\tilde{u}\|_{L^2}^2
   +
   \frac{\eta}{8}\|\Delta\tilde{b}\|_{L^2}^2,
\end{align*}
and
\begin{align*}
   |\Lambda_{21} | 
   &
   \leq
   \|\tilde{b}\|_{L^6}\|\nabla u'\|_{L^3}\|\Delta\tilde{b}\|_{L^2}
   \leq
   C\|\Delta u'\|_{L^2}^2\|\nabla\tilde{b}\|_{L^2}^2
   +
   \frac{\eta}{16}\|\Delta\tilde{b}\|_{L^2}^2.
\end{align*}
Combining all the above estimates, and by denoting 
\begin{align*}
  E_1(t) 
  &
  = 
  \|\tilde{u}(t)\|_{L^2}^2 
  + 
  \alpha^2\|\nabla\tilde{u}(t)\|_{L^2}^2 
  +
  \|\nabla\tilde{u}(t)\|_{L^2}^2
  \\&\quad 
  +
  \|\tilde{w}(t)\|_{L^2}^2 
  + 
  \|\tilde{b}(t)\|_{L^2}^2 
  + 
  \|\nabla \tilde{b}(t)\|_{L^2}^2,
\end{align*}
and
\begin{align*}
  E_2(t) 
  &
  = 
  \nu\|\nabla\tilde{u}(t)\|_{L^2}^2
  +
  \nu\|\Delta\tilde{u}(t)\|_{L^2}^2
  +
  \nu\|\nabla\tilde{w}(t)\|_{L^2}^2
  \\&\quad
  +
  \eta\|\nabla\tilde{b}(t)\|_{L^2}^2
  +
  \eta\|\Delta\tilde{b}(t)\|_{L^2}^2,
\end{align*}
we obtain 
\begin{align*}
  \frac{d E_{1}(t)}{dt} 
  +
  E_{2}(t)
  \leq
  ME_{1}(t),
\end{align*}
where the constant $M$ depends on $\nu, \eta$, and the $L^2$-, $H^1$-, and $H^2$-norms of $u, u', b, b'$, and the $L^2$- and $H^1$-norms of $w, w'$. 

Therefore, by the initial condition $\tilde{u}_0=\tilde{w}_0=\tilde{b}_0=0$ (thus $\nabla\tilde{u}_0=0$ and $\nabla\tilde{b}_0=0$), Gr\"onwall inequality implies that $\tilde{u}=\tilde{w}=\tilde{b}=0$, i.e., 
$$u=u', w=w', b=b'.$$
Thus, the proof of uniqueness is now complete.  
\qedsymbol

\section{Proof of the Convergence Theorems}\label{conv}

\subsection{Proof of Theorems~\ref{thm-conv1}}\label{proof-conv1}
We begin by taking the {\it curl} operator ``$\nabla\times$'' to the $u$-equation in System~\eqref{VVV_MHD}, denoting by $\omega=\nabla\times u$, and obtain 
\begin{align}
   &(I - \alpha^2\Delta)\frac{\partial \omega}{\partial t} 
   - 
   \nu\Delta \omega 
   + 
   (u\cdot\nabla)w
   - 
   (w\cdot\nabla)u
   = 
   \nabla\times(j\times b),
\end{align}
where we used the divergene-free condition of $u$ and $w$. 
Then, we denote by $\xi=\omega-w$, and subtract the $w$-equation in \eqref{MHD} from the above equation of $\omega$, and have
\begin{equation*}
   (I - \alpha^2\Delta)\frac{\partial \xi}{\partial t}
   - 
   \nu\Delta \xi
   =
  \alpha^2\Delta\frac{\partial w}{\partial t}.
\end{equation*}
Multiplying the above equation of $\xi$ by $\xi$ and integrating by parts over $\mathbb{T}^3$, we obtain 
\begin{align*}
  &
  \frac{1}{2}\frac{d}{dt}\big(\|\xi\|_{L^2}^2 + \alpha^2\|\nabla\xi\|_{L^2}^2\big)
  +
  \|\nabla\xi\|_{L^2}^2
  = 
  \alpha^2\int_{\mathbb{T}^3}\Delta \frac{\partial w}{\partial t}\cdot\xi
  =
  \alpha^2\int_{\mathbb{T}^3}\frac{\partial w}{\partial t}\cdot\Delta\xi
  \\&
  =
  \alpha^2\nu\int_{\mathbb{T}^3}\Delta w\cdot\Delta\xi
  -
  \alpha^2\int_{\mathbb{T}^3}(u\cdot\nabla)w\cdot\Delta\xi 
  \\&\quad
  +
  \alpha^2\int_{\mathbb{T}^3}(w\cdot\nabla)u\cdot\Delta\xi 
  +
  \alpha^2\int_{\mathbb{T}^3}\nabla\times(j\times b)\cdot\Delta\xi,
\end{align*}
where we substituted the $w$-equation of \eqref{MHD} for the time-derivative of $w$. 

Next, we estimate the four integrals on the right side of the above equation. For the first one, we integrate by parts again and use Cauchy-Schwarz inequality, and have
\begin{align*}
  &
  \Big|\alpha^2\nu\int_{\mathbb{T}^3}\Delta w\cdot\Delta\xi\Big|
  =
  \Big|\alpha^2\nu\int_{\mathbb{T}^3}\nabla\Delta w\cdot\nabla\xi\Big|
  \\&
  \leq
  C\|w\|_{H^3}\|\nabla\xi\|_{L^2}
  \leq
  C\alpha^2\nu^2\|w\|_{H^3}^2
  +
  \alpha^2\|\nabla\xi\|_{L^2}^2.
\end{align*}
For the second integral, we also integrate by parts and apply H\"older and Sobolev inequalities, and obtain
\begin{align*}
  &
  \Big|-\alpha^2\int_{\mathbb{T}^3}(u\cdot\nabla)w\cdot\Delta\xi\Big|
  \leq
  \alpha^2\int_{\mathbb{T}^3}|\nabla u||w||\nabla\xi|
  +
  \alpha^2\int_{\mathbb{T}^3}|u||\Delta w||\nabla\xi|
  \\&
  \leq
  C\alpha^2\|\nabla u\|_{L^3}\|\nabla w\|_{L^6}\|\nabla\xi\|_{L^2}
  +
  C\alpha^2\|u\|_{L^{\infty}}\|\Delta w\|_{L^2}\|\nabla\xi\|_{L^2}
  \\&
  \leq
  C\alpha^2\|u\|_{H^2}\|w\|_{H^2}\|\nabla\xi\|_{L^2}
  \leq
  C\alpha^2\|u\|_{H^2}^2\|w\|_{H^2}^2
  +
  \alpha^2\|\nabla\xi\|_{L^2}^2,
\end{align*}
where we also used Agmon's inequality \eqref{agmon1}. By similar arguments, the third integral is also bounded by 
\begin{align*}
  &
  \Big|\alpha^2\int_{\mathbb{T}^3}(w\cdot\nabla)u\cdot\Delta\xi\Big|
  \leq
  C\alpha^2\|u\|_{H^2}^2\|w\|_{H^2}^2
  +
  \alpha^2\|\nabla\xi\|_{L^2}^2.
\end{align*}
As for the last integral, we integrate by parts again and obtain
\begin{align*}
  &
  \Big|\alpha^2\int_{\mathbb{T}^3}\nabla\times(j\times b)\cdot\Delta\xi,\Big|
  \leq
  \alpha^2\int_{\mathbb{T}^3}|\Delta(j\times b)||\nabla\xi|
  \\&
  \leq
  C\alpha^2\int_{\mathbb{T}^3}|\Delta j||b||\nabla\xi|
  +
  C\alpha^2\int_{\mathbb{T}^3}|j||\Delta b||\nabla\xi|
  +
  C\alpha^2\int_{\mathbb{T}^3}|\nabla j||\nabla b||\nabla\xi|
  \\&
  \leq
  C\alpha^2\|\Delta j\|_{L^2}\|b\|_{L^{\infty}}\|\nabla\xi\|_{L^2}
  +
  C\alpha^2\|j\|_{L^3}\|\Delta b\|_{L^6}\|\nabla\xi\|_{L^2}
  \\&\quad
  +
  C\alpha^2\|\nabla j\|_{L^3}\|\nabla b\|_{L^6}\|\nabla\xi\|_{L^2}
  \\&
  \leq
  C\alpha^2\|b\|_{H^3}^4
  +
  \alpha^2\|\nabla\xi\|_{L^2}^2,
\end{align*}
where we used H\"older and Sobolev inequalities as well as Cauchy-Schwarz inequality. 

Finally, by combining all the above estimates, we have 
\begin{align*}
  &
  \frac{d}{dt}\big(\|\xi\|_{L^2}^2 + \alpha^2\|\nabla\xi\|_{L^2}^2\big)
  \leq
  C\big(\|\xi\|_{L^2}^2 + \alpha^2\|\nabla\xi\|_{L^2}^2\big)
  +
  K\alpha^2,
\end{align*}
where the constant $K_0$ depends on the constant $C$ in the estimates above, as well as on the $H^2$-norm of $u$ and the $H^3$-norms of $w$ and $b$. Therefore, Gr\"onwall inequality implies 
\begin{align*}
  &
  \|\xi(t)\|_{L^2}^2
  +
  \alpha^2\|\nabla\xi(t)\|_{L^2}^2
  +
  \int_{0}^{t}\|\nabla\xi(s)\|_{L^2}^2\,ds
  \\&
  \leq
  K_0\alpha^2\big(e^{Ct} - 1\big).
\end{align*}
Since $\xi(\cdot, 0)=\xi_0=0$, we have 
\begin{align*}
  &
  \|\omega(t) - w(t)\|_{L^{\infty}(0, T; H)}^2 + \|\nabla\omega(t) - \nabla w(t)\|_{L^2(0,T;V)}^2
  \leq 
  C_{K_0}\alpha^2e^{CT}. 
\end{align*}
Proof of our first convergence result is now complete. 
\qedsymbol

\subsection{Proof of Theorems~\ref{thm-conv2}}\label{proof-conv2}

Note that, by Theorem~\ref{thm-conv1}, it suffices to prove $w\to\nabla\times U:=\Omega$ in $L^2$ as $\alpha\to 0$. In order to do so, and to prove other convergence in Theorem~\ref{thm-conv2}, we first take the curl operator ``\,$\nabla\times$\,'' to the $U$-equation in System~\eqref{MHD}, and get
\begin{align*}
  &
  \frac{\partial \Omega}{\partial t} 
  - 
  \nu\Delta \Omega 
  + 
  (U\cdot\nabla)\Omega
  -
  (\Omega\cdot\nabla)U
  = 
  \nabla\times(J\times B),
\end{align*}
where $J=\nabla\times B$ is the density of the induced electric current. 
Then, denote by $\zeta=u-U$, $q=w-\Omega$, $\beta=b-B$, $\mu=j-J=\nabla\times\beta$, and $\tilde{\Pi}=p-P$, and subtract the above equation from the $w$-equation in System~\eqref{VVV_MHD}, as well as subtract the $u$- and $B$-equations in System~\eqref{MHD} from the $u$- and $b$-equations in System~\eqref{VVV_MHD}, respectively, and subtract the curl of the $u$- and $B$-equations in System~\eqref{MHD} from the curl of the $u$- and $b$-equations in System~\eqref{VVV_MHD}, respectively, and we obtain
\begin{equation}
 \left\{
 \begin{aligned}
  &
  (I - \alpha^2\Delta)\frac{\partial \zeta}{\partial t} 
  -
  \alpha^2\Delta \frac{\partial U}{\partial t}
  - 
  \nu\Delta \zeta
  + 
  w\times\zeta
  +
  q\times U
  +
  \nabla \tilde{\Pi} 
  \\&\quad
  = 
  \mu\times b
  +
  J\times\beta,
  \\&
  \frac{\partial q}{\partial t}
  -
  \nu\Delta q
  +
  (\zeta\cdot\nabla)w
  +
  (U\cdot\nabla)q
  -
  (w\cdot\nabla)\zeta
  -
  (q\cdot\nabla)U
  \\&\quad
  =
  (b\cdot\nabla)\mu
  +
  (\beta\cdot\nabla)J
  -
  (j\cdot\nabla)\beta
  -
  (\mu\cdot\nabla)B,
  \\&
  \frac{\partial \beta}{\partial t}
  -
  \eta\Delta \beta
  +
  (\zeta\cdot\nabla)b
  +
  (U\cdot\nabla)\beta
  =
  (b\cdot\nabla)\zeta
  +
  (\beta\cdot\nabla)U,
  \\&
  \frac{\partial \mu}{\partial t}
  -
  \eta\Delta\mu
  +
  \nabla\times\big((\zeta\cdot\nabla)b + (U\cdot\nabla)\beta\big)
  =
  \nabla\times\big((b\cdot\nabla)\zeta + (\beta\cdot\nabla)U\big),
  \\&
  0=\nabla\cdot\zeta=\nabla\cdot\beta=\nabla\cdot q,
  \\&
  \zeta(\cdot,0)=\zeta_0=0, \quad q(\cdot,0)=q_0=0,\quad \beta(\cdot,0)=\beta_0=0.
 \end{aligned}
 \right.
 \label{Sys3}
\end{equation} 
Next, we multiply the $\zeta$-, $q$-, $\beta$-, and $\mu$-equations in System~\eqref{Sys3} by $\zeta$, $q$, $\beta$, and $\mu$, respectively, integrate over $\mathbb{T}^3$, and add, so that we have
\begin{align*}
  &
  \frac{1}{2}\frac{d}{dt}\big(\|\zeta\|_{L^2}^2 + \alpha^2\|\nabla\zeta\|_{L^2}^2 + \|q\|_{L^2}^2 + \|\beta\|_{L^2}^2 + \|\mu\|_{L^2}^2\big)
  \\&\quad
  +
  \nu\big(\|\nabla\zeta\|_{L^2}^2 + \|\nabla q\|_{L^2}^2\big)
  +
  \eta\big(\|\nabla\beta\|_{L^2}^2 + \|\nabla\mu\|_{L^2}^2\big)
  \\&
  =
  \alpha^2\int \Delta\frac{\partial U}{\partial t}\cdot\zeta
  -
  \int (q\times U)\cdot\zeta
  +
  \int (\mu\times b)\cdot\zeta
  +
  \int (J\times \beta)\cdot\zeta
  \\&\quad
  -
  \int (\zeta\cdot\nabla)w\cdot q
  +
  \int (w\cdot\nabla)\zeta\cdot q
  +
  \int (q\cdot\nabla)U\cdot q
  \\&\quad\quad
  +
  \int (b\cdot\nabla)\mu\cdot q
  +
  \int (\beta\cdot\nabla)J\cdot q
  -
  \int (j\cdot\nabla)\beta\cdot q
  -
  \int (\mu\cdot\nabla)B\cdot q
  \\&\quad
  -
  \int (\zeta\cdot\nabla)b\cdot\beta
  +
  \int (b\cdot\nabla)\zeta\cdot\beta
  +
  \int (\beta\cdot\nabla)U\cdot\beta
  \\&\quad
  +
  \int (\zeta\cdot\nabla)b\cdot(\nabla\times\mu)
  +
  \int (U\cdot\nabla)\beta\cdot(\nabla\times\mu)
  \\&\quad\quad
  -
  \int (b\cdot\nabla)\zeta\cdot(\nabla\times\mu)
  -
  \int (\beta\cdot\nabla)U\cdot(\nabla\times\mu)
  \\&
  :=
  \Theta_1 + \Theta_2 + \cdots + \Theta_{18},
\end{align*}
where we also integrated by part on the right side of the above equation. 

Then, we estimate the eighteen integrals on the right side of the above equation. For the first one, we substitute $\cfrac{\partial U}{\partial t}$ by the $U$-equation in System~\eqref{MHD}, integrate by parts, and obtain 
\begin{align*}
  |\Theta_1|
  &
  \leq
  \alpha^2\nu\int |\nabla\Delta U||\nabla\zeta|
  +
  \alpha^2\int |\nabla(\Omega\times U)||\nabla\zeta|
  +
  \alpha^2\int |\nabla(J\times B)||\nabla\zeta|
  \\&
  \leq
  C\alpha^2\nu\|U\|_{H^3}^2
  +
  C\alpha^2\|U\|_{H^2}^2
  +
  C\alpha^2\|U\|_{H^1}\|U\|_{H^3}
  \\&\quad
  +
  C\alpha^2\|B\|_{H^2}^2
  +
  C\alpha^2\|B\|_{H^1}\|B\|_{H^3}
  +
  \frac{\nu}{16}\|\nabla\zeta\|_{L^2}^2
  \\&
  :=\alpha^2 K_1 + \frac{\nu}{16}\|\nabla\zeta\|_{L^2}^2. 
\end{align*}
In order to estimates $\Theta_2$, we apply H\"older and Sobolev inequalities, and get
\begin{align*}
  |\Theta_2|
  &
  \leq
  C\|q\|_{L^3}\|U\|_{L^6}\|\zeta\|_{L^2}
  \leq
  C\|\nabla q\|_{L^2}\|\nabla U\|_{L^2}\|\zeta\|_{L^2}
  \\&
  \leq
  C\|\nabla U\|_{L^2}^2\|\zeta\|_{L^2}^2 + \frac{\nu}{16}\|\nabla q\|_{L^2}^2,
\end{align*}
where we also used Poincar\'e's and Young's inequalities. 

Estimates of $\Theta_3$ is similar, and we have 
\begin{align*}
  |\Theta_3|
  \leq
  C\|\mu\|_{L^3}\|b\|_{L^2}\|\zeta\|_{L^6}
  \leq
  C\|b\|_{L^2}^2\|\nabla\zeta\|_{L^2}^2 + \frac{\eta}{16}\|\nabla\mu\|_{L^2}^2,
\end{align*}
and 
\begin{align*}
  |\Theta_4|
  \leq
  C\|J\|_{L^2}^2\|\mu\|_{L^2}^2 + \frac{\nu}{16}\|\nabla\zeta\|_{L^2}^2,
\end{align*}
where we used the fact that $\mu=\nabla\times\beta$ thus $\|\nabla\beta\|_{L^2}\leq C\|\mu\|_{L^2}$. 

As for $\Theta_5$, we substitute $w=q+\Omega$ and take advantage of the symmetry of the nonlinear term which leads to the cancellation of one of the two integrals, i.e. we obtain, 
\begin{align*}
  |\Theta_5|
  &
  =
  \Big|\int (\zeta\cdot\nabla)(q+\Omega)\cdot q\Big|
  =
  \Big|\int (\zeta\cdot\nabla)\Omega\cdot q\Big|
  \\&
  \leq
  C\|\zeta\|_{L^6}\|\nabla\Omega\|_{L^3}\|q\|_{L^2}
  \leq
  C\|\Delta \Omega\|_{L^2}^2\|q\|_{L^2}^2
  + 
  \frac{\nu}{16}\|\nabla\zeta\|_{L^2}^2.
\end{align*}

Then, we proceed to estimates $\Theta_6$ as
\begin{align*}
  |\Theta_6|
  \leq
  C\|w\|_{L^6}\|\nabla\zeta\|_{L^2}\|q\|_{L^3}
  \leq
  C\|\nabla w\|_{L^2}^2\|\nabla\zeta\|_{L^2}^2 + \frac{\nu}{16}\|\nabla q\|_{L^2}^2,
\end{align*}
while for $\Theta_7$, we apply H\"older and Sobolev inequalities, and get
\begin{align*}
  |\Theta_7|
  \leq
  C\|q\|_{L^{2}}\|\nabla U\|_{L^3}\|q\|_{L^6}
  \leq
  C\|\Delta U\|_{L^2}^2\|q\|_{L^2}^2 + \frac{\nu}{16}\|\nabla q\|_{L^2}^2.
\end{align*} 

The next integral $\Theta_8$ is bounded by
\begin{align*}
  |\Theta_8|
  &
  \leq
  C\|b\|_{L^{6}}\|\nabla\mu\|_{L^{2}}\| q\|_{L^3}
  \leq
  C\|\nabla b\|_{L^{2}}\|\nabla\mu\|_{L^2}\|q\|_{L^2}^{1/2}\|\nabla q\|_{L^2}^{1/2}
  \\&
  \leq
  C\|\nabla b\|_{L^{2}}^4\|q\|_{L^2}^2
  + \frac{\nu}{16}\|\nabla q\|_{L^2}^2
  + \frac{\eta}{16}\|\nabla\mu\|_{L^2}^2,
\end{align*}
and we also point out that the $L^2$ in time integral of $\|\nabla b\|_{L^2}$ is independent of $\alpha$, since the $L^2$ in time integral of $\|\nabla u\|_{L^2}$ is independent of $\alpha$. 

Regarding $\Theta_9$, we have 
\begin{align*}
  |\Theta_9|
  \leq
  C\|\beta\|_{L^{3}}\|\nabla J\|_{L^{2}}\| q\|_{L^6}
  \leq
  C\|\nabla J\|_{L^2}^4\|\beta\|_{L^2}^2 
  + \frac{\nu}{16}\|\nabla q\|_{L^2}^2
  + \frac{\eta}{16}\|\nabla\beta\|_{L^2}^2.
\end{align*}
 
As for $\Theta_{10}$, we use $\|\nabla\beta\|_{L^2}\leq C\|\mu\|_{L^2}$, and get
\begin{align*}
  |\Theta_{10}|
  &
  \leq
  C\|j\|_{L^{2}}\|\nabla\beta\|_{L^{6}}\| q\|_{L^3}
  \leq
  \|\nabla b\|_{L^2}\|\nabla\mu\|_{L^2}\|q\|_{L^2}^{1/2}\|\nabla q\|_{L^2}^{1/2}
  \\&
  \leq
  C\|\nabla b\|_{L^2}^4\|q\|_{L^2}^2 
  + \frac{\nu}{16}\|\nabla q\|_{L^2}^2
  + \frac{\eta}{16}\|\nabla\mu\|_{L^2}^2,
\end{align*}
where we note that $\|\nabla b\|_{L^2}$ is $L^2$-integrable in time. 

Estimates of $\Theta_{11}$ is similar and we have
\begin{align*}
  |\Theta_{11}|
  \leq
  C\|\mu\|_{L^6}\|\nabla B\|_{L^{2}}\|q\|_{L^3}
  \leq
  C\|\nabla B\|_{L^2}^4\|q\|_{L^2}^2 
  + \frac{\nu}{16}\|\nabla q\|_{L^2}^2
  + \frac{\eta}{16}\|\nabla\mu\|_{L^2}^2.
\end{align*}
Estimates for the next two terms are similar. For $\Theta_{12}$, we first integrate by parts, then estimate as  
\begin{align*}
  |\Theta_{12}|
  &
  \leq
  \int |\zeta||b||\nabla\beta|
  \leq
  C\|\zeta\|_{L^6}\|b\|_{L^{2}}\|\mu\|_{L^3}
  \\&
  \leq
  C\|b\|_{L^2}^4\|\mu\|_{L^2}^2 
  + \frac{\nu}{16}\|\nabla\zeta\|_{L^2}^2
  + \frac{\eta}{16}\|\nabla\mu\|_{L^2}^2,
\end{align*}
while for $\Theta_{13}$, we have
\begin{align*}
  |\Theta_{13}|
  &
  \leq
  \int |b||\zeta||\nabla\beta|
  \leq
  C\|\zeta\|_{L^6}\|b\|_{L^{2}}\|\mu\|_{L^3}
  \\&
  \leq
  C\|b\|_{L^2}^4\|\mu\|_{L^2}^2 
  + \frac{\nu}{16}\|\nabla\zeta\|_{L^2}^2
  + \frac{\eta}{16}\|\nabla\mu\|_{L^2}^2,
\end{align*}
Regarding $\Theta_{14}$, we estimate as
\begin{align*}
  |\Theta_{14}|
  \leq
  C\|\nabla U\|_{L^2}^4\|\beta\|_{L^2}^2 + \frac{\eta}{16}\|\nabla\beta\|_{L^2}^2.
\end{align*}
and
\begin{align*}
  |\Theta_{15}|
  &
  \leq
  C\|\zeta\|_{L^{6}}\|\nabla b\|_{L^{3}}\|\nabla\mu\|_{L^2}
  \\&
  \leq
  C\|\nabla b\|_{L^2}\|\Delta b\|_{L^2}\|\zeta\|_{L^2}^2 
  + \frac{\eta}{16}\|\nabla\mu\|_{L^2}^2. 
\end{align*}
In order to estimate $\Theta_{16}$, we apply Agmon's inequality~\eqref{agmon1}, and get
\begin{align*}
  |\Theta_{16}|
  \leq
  C\|U\|_{L^{\infty}}\|\nabla\beta\|_{L^{2}}\|\nabla\mu\|_{L^2}
  \leq
  C\|U\|_{H^2}^2\|\mu\|_{L^2}^2 
  + \frac{\eta}{16}\|\nabla\mu\|_{L^2}^2,
\end{align*}
where we also used the fact $\|\nabla\beta\|_{L^2}\leq C\|\mu\|_{L^2}$. 

Estimates of $\Theta_{17}$ is similar to $\Theta_{16}$, so we have
\begin{align*}
  |\Theta_{17}|
  \leq
  C\|b\|_{L^{\infty}}\|\nabla\zeta\|_{L^{2}}\|\nabla\mu\|_{L^2}
  \leq
  C\|\Delta b\|_{L^2}^2\|\nabla\zeta\|_{L^2}^2 
  + \frac{\eta}{16}\|\nabla\mu\|_{L^2}^2.
\end{align*}
Finally, the last term $\Theta_{18}$ is bounded by 
\begin{align*}
  |\Theta_{18}|
  \leq
  C\|\beta\|_{L^{6}}\|\nabla U\|_{L^{3}}\|\nabla\mu\|_{L^2}
  \leq
  C\|\Delta U\|_{L^2}^2\|\mu\|_{L^2}^2 
  + \frac{\eta}{16}\|\nabla\mu\|_{L^2}^2. 
\end{align*}

Combining all the above estimates, after some cancellation and rearrangements, and by denoting 
$$E(t):=\big(\|\zeta\|_{L^2}^2 + \alpha^2\|\nabla\zeta\|_{L^2}^2 + \|q\|_{L^2}^2 + \|\beta\|_{L^2}^2 + \|\mu\|_{L^2}^2\big)^{1/2}$$
we obtain 
\begin{align*}
  &
  \frac{1}{2}\frac{d}{dt}E(t)^2
  +
  \nu\big(\|\nabla\zeta\|_{L^2}^2 + \|\nabla q\|_{L^2}^2\big)
  +
  \eta\big(\|\nabla\beta\|_{L^2}^2 + \|\nabla\mu\|_{L^2}^2\big)
  \\&
  \leq
  \alpha^2 K_1
  +
  K_2\big(\|b\|_{H^1}^2+\|b\|_{H^2}^2\big)E(t)^2.
\end{align*}
Thanks to the $L^2$ in time integrability of the $H^1$- and $H^2$-norms of $b$, Gr\"onwall inequality implies that 
\begin{align*}
  &
  \|\zeta(T)\|_{L^2}^2 + \alpha^2\|\nabla\zeta(T)\|_{L^2}^2 + \|q(T)\|_{L^2}^2 + \|\beta(T)\|_{L^2}^2 + \|\mu(T)\|_{L^2}^2
  \\&\quad
  +
  \nu\int_{0}^{T} \big(\|\nabla\zeta(t)\|_{L^2}^2 + \|\nabla q(t)\|_{L^2}^2\big)\,dt
  \\&\quad\quad
  +
  \eta\int_{0}^{T} \big(\|\nabla\beta(t)\|_{L^2}^2 + \|\nabla\mu(t)\|_{L^2}^2\big)\,dt
  \\&
  \leq
  C_{T}\alpha^2.
\end{align*}
Proof of our second convergence result is now complete. 
\qedsymbol

\subsection{Proof of Theorem~\ref{thm-blowup}}\label{proof-blowup}
We argue by contradiction, so suppose that strong solutions $U, B\in L^{\infty}(0, T; H^3)\cap L^2(0, T; H^4)$ to system \eqref{MHD} exist. In the energy identity (\eqref{Energy}), we take $\limsup$ as $\alpha\to 0^{+}$, which gives 
    \begin{align*}
        &
        \alpha^2\normL{\nabla U(t)}^2
        +
        \normL{U(t)}^2
        +
        \normL{B(t)}^2
        +
        2\int_{0}^{t}\big(\nu\normL{\nabla U(s)}^2 + \eta\normL{\nabla B(s)}^2\big)\,ds
        \\& \nonumber
        =
        \normL{U_0}^2
        +
        \normL{B_0}^2.
    \end{align*}
However, the energy identity of system \eqref{MHD} states that 
    \begin{align*}
        &
        \normL{U(t)}^2
        +
        \normL{B(t)}^2
        +
        2\int_{0}^{t}\big(\nu\normL{\nabla U(s)}^2 + \eta\normL{\nabla B(s)}^2\big)\,ds
        \\& \nonumber
        =
        \normL{U_0}^2
        +
        \normL{B_0}^2, 
    \end{align*}
which leads to a contradiction, due to the hypothesis (\eqref{blow-up_condition}). 
\qedsymbol

\section*{Declarations}

\subsection*{Data Availability Statement}
No data was generated in this study.

\subsection*{Conflict of Interest Statement}
On behalf of authors A.~L. and Y.~P., the corresponding author states that there is no conflict of interest. 

\subsection*{Author Contribution Statement}
Authors A.~L. and Y.~P. contributed equally to the paper. 

\subsection*{Funding Declaration}
The research of A.L. was supported in part by NSF grants DMS-2206762, DMS-2510494, CMMI-1953346, and USGS Grant No. G23AC00156-01.

\end{document}